\documentclass[11pt]{article}
\usepackage{amsmath}
\usepackage{amsmath,amsthm,amsfonts,amssymb,bm,euscript}
\usepackage{fancyhdr,amscd}
\usepackage{anysize}
\usepackage{graphicx,epsfig}
\usepackage{amsthm,latexsym,amsmath,amssymb}
\usepackage{amssymb,amsmath}
\usepackage{mathrsfs}
\usepackage{graphicx}
\usepackage{color}
\usepackage{mathrsfs}
\usepackage{cite}
\usepackage{indentfirst}
 \usepackage[titletoc]{appendix}
 \usepackage [latin1]{inputenc}
\let\d=\delta

\let\D=\Delta

\def\cC{{\cal C}}

\def\eqdefa{\buildrel\hbox{\footnotesize def}\over =}

\def\C{\mathop{\bf C\kern 0pt}\nolimits}
\def\DD{\mathop{\bf D\kern 0pt}\nolimits}
\def\K{\mathop{\bf K\kern 0pt}\nolimits}
\def\N{\mathop{\bf N\kern 0pt}\nolimits}
\def\Q{\mathop{\bf Q\kern 0pt}\nolimits}
\def\R{\mathop{\bf R\kern 0pt}\nolimits}

\renewcommand{\div}{\mbox{\rm div}\;\!}

\def\cC{{\cal C}}

\def\lg{\mbox{log}}

\def\eqdefa{\buildrel\hbox{\footnotesize def}\over =}

\def\C{\mathop{\bf C\kern 0pt}\nolimits}
\def\DD{\mathop{\bf D\kern 0pt}\nolimits}
\def\K{\mathop{\bf K\kern 0pt}\nolimits}
\def\N{\mathop{\bf N\kern 0pt}\nolimits}
\def\Q{\mathop{\bf Q\kern 0pt}\nolimits}
\def\R{\mathop{\bf R\kern 0pt}\nolimits}

\renewcommand{\div}{\mbox{\rm div}\;\!}

\newcommand{\Cl}{{\rm curl}}
\newcommand{\Dv}{{\rm div}}

\newcommand{\beq}{\begin{equation}}
\newcommand{\eeq}{\end{equation}}
\newcommand{\ben}{\begin{eqnarray}}
\newcommand{\een}{\end{eqnarray}}
\newcommand{\beno}{\begin{eqnarray*}}
\newcommand{\eeno}{\end{eqnarray*}}

\newtheorem{Theorem}{Theorem}[section]

\newtheorem{Proposition}[Theorem]{Proposition}
\newtheorem{Lemma}[Theorem]{Lemma}

\newtheorem{Remark}[Theorem]{Remark}

\numberwithin{equation}{section}

\allowdisplaybreaks \numberwithin{equation} {section}

\begin{document}
\title{The unique global solvability of the nonhomogeneous incompressible  asymmetric fluids with vacuum  \thanks {Research supported by the
National Natural Science Foundation of China (11501332,11771043,51976112), the  Natural Science Foundation of Shandong Province (ZR2015AL007),
 and Young Scholars Research Fund of Shandong University of Technology.}
}
\author{ Fuyi  Xu$^\dag$  \ Mingxue  Zhang \ Liening  Qiao \\[2mm]
 { \small  School of Mathematics and  Statistics, Shandong University of Technology,}\\
  { \small Zibo,    255049,  Shandong,    China}}
         \date{}
         \maketitle
\noindent{\bf Abstract}\ \ \ The present paper deals with
the nonhomogeneous incompressible  asymmetric fluids equations in dimension $d= 2,3$. The aim is to prove the unique global solvability  of the system with only bounded nonnegative initial density and  $H^{1}$ initial velocities. We first construct the global existence of the solution with large data in 2-D. Next, we establish the existence of local in time solution for arbitrary large data and global in time  for some smallness conditions in 3-D. Finally,  the uniqueness of the solution is proved under quite soft assumptions about its regularity through a Lagrangian approach. In particular, the initial vacuum is allowed.
\vskip   0.2cm \footnotetext[1]{$^\dag$Corresponding author.}
\vskip   0.2cm \footnotetext[2]{E-mail addresses: zbxufuyi@163.com(F.Xu)\ .} \setlength{\baselineskip}{20pt}

\section{Introduction and Main Results}
\setcounter{section}{1}\setcounter{equation}{0} \ \ \ \ \
 \par
In the present paper, we consider the following $d$-dimensional (for $d= 2,3$) nonhomogeneous incompressible  asymmetric fluids equations in an open bounded set $\Omega$:
\begin{align}\label{eq:NA}
\left\{
\begin{aligned}
&\rho(u_{t}+u\cdot\nabla u)-(\mu+\chi)\Delta u+\nabla P=2\chi\Cl\omega &\quad \text{in}\quad  \mathbb{R}^{+}\times \Omega,\\
&\Dv u=0 &\quad \text{in}\quad  \mathbb{R}^{+}\times \Omega, \\
&\rho(\omega_{t}+u\cdot\nabla \omega)-\gamma\Delta\omega-\kappa\nabla\Dv\omega+4\chi\omega
=2\chi \Cl u &\quad \text{in}\quad  \mathbb{R}^{+}\times \Omega,\\
&\rho_{t}+u\cdot\nabla\rho=0 &\quad \text{in}\quad  \mathbb{R}^{+}\times \Omega,\\
&(\rho,u,\omega)|_{t=0}=(\rho_{0},u_{0},\omega_{0}) &\text{in}\quad  \Omega,
\end{aligned}
\right.
\end{align}
where $u$ is  the fluid velocity, $\omega$ is the field of microrotation representing the angular velocity of the rotation of the particles of the fluid, $P$ is the scalar pressure of the flow,  while $\rho_0$, $u_0$ and $\omega_0$ are the given initial density, initial velocity and initial angular velocity, respectively, with $\textrm{div} u_0=0$. $\mu$ is the Newtonian kinematic viscosity,  $\kappa$ is  the angular viscosity, $\chi$ is the micro-rotation viscosity. For the derivation of   system \eqref{eq:NA} and a discussion on their physical meaning, see \cite{CD}. Concerning applications, the
micropolar fluid model has been used, for example, in lubrication theory \cite{CF,PS}, as well as in modeling blood flow in thin
vessels \cite{AT}.  We shall  assume that the fluid domain $\Omega$ is either
the   torus $\mathbb{T}^d$ or a $\cC^2$ simply connected bounded domain of $\mathbb{R}^d$.
For simplicity, we take $\nu=\mu+\chi$.

System \eqref{eq:NA} includes several important models as special cases.   When $\rho=const$,  system \eqref{eq:NA} becomes  the  incompressible micropolar fluid provided that $P$ is an unknown pressure function, which was first introduced in 1965 by Eringen to model micropolar fluids (see Eringen \cite{ACE}, Sections 1 and 6). The micropolar fluid model that will be considered in this article is a generalization of the Navier-Stokes equations,  which takes into account the microstructure of the fluid by which we mean the geometry and microrotation of particles.  As experiments show the solutions of this model represent flows of many fluids (like,
e.g. blood, see Reference \cite{PRU}) more precisely than solutions of the Navier-Stokes equations. Due to its physical significance and mathematical relevance, there have a lot of works studying on the viscous or inviscid  3-D(2-D) system (see e.g. \cite{CM,DZ,DLW,XUE}). When $\omega=0$, system \eqref{eq:NA} reduces the  nonhomogeneous incompressible   Navier-Stokes equations,  which is obtained by mixing two miscible fluids that
are incompressible and that have different densities. It may also describe a fluid
containing a melted substance. One may check \cite{PLS} for the detailed derivation of
this system. Kazhikov \cite{Kaz} first proved the global existence of weak solutions for the  3-D nonhomogeneous incompressible  Navier-Stokes equations with $\inf\rho_{0}>0$.  Later on, Simon \cite{Sim2} removed the lower bound assumption on $\rho_{0}$, and Lions \cite{PLS} proved that $\rho$ is a renormalized solution of the mass equation. However, the uniqueness and smoothness of weak solutions to the nonhomogeneous incompressible Navier-Stokes equations, even for the 2-D case, is still an open problem. Local existence (but without uniqueness) of strong solutions to the nonhomogeneous
incompressible Navier-Stokes equations was first established by Antontsev and Kazhikov \cite{AK}, under the assumptions that the initial density is bounded and away
from zero and the initial velocity has $H^{1}$ regularity. In the bounded domain $\Omega$, Ladyzhenskaya and Solonnikov \cite{LS} first constructed global strong solutions in 2-D,  and unique local in time maximal strong solution for arbitrary data and global in time for small data in  3-D  where $\rho_{0}\in C^{1}(\Omega)$ is bounded away from zero. In general when $\rho_{0}\in L^{\infty}(\mathbb{R}^{d})$ with a positive lower bound and $u_{0}\in H^{2}(\mathbb{R}^{d})$, Danchin and Mucha \cite{Dan3} proved that the nonhomogeneous incompressible  Navier-Stokes equations have a unique local strong solution. Furthermore, with the initial density fluctuation being sufficiently small, they also obtained the global well-posedness.  Paicu et al. \cite{PZZ} proved the global existence and uniqueness of the solution to
$d$-dimensional (for $d = 2,3$) nonhomogeneous incompressible Navier-Stokes equations
with initial density being bounded from above and below by some positive constants,
and with initial velocity $u_0\in H^{s}(\mathbb{R}^{2})(s>0)$ in 2-D, or $u_0\in H^{1}(\mathbb{R}^{3})$ satisfying some smallness conditions in 3-D.  If it is not assumed that the density is bounded away from zero, then the analysis of nonhomogeneous incompressible Navier-Stokes equations
gets wilder, since the initial density is allowed to have a vacuum. Choe-Kim
\cite{CHK} first proved the local existence and uniqueness of strong solution with initial data $(\rho_0,u_0)$ satisfying $0\leq\rho_0\in L^{\infty}(\mathbb{R}^{3})\cap H^{1}(\mathbb{R}^{3}),  u_0\in H^{2}(\mathbb{R}^{3})\cap H^{1}_{0,\sigma}(\mathbb{R}^{3})$  and the compatibility conditions.  Recently, in  a bounded domain $\Omega$,  removed the compatibility condition,  Li \cite{Lijinkai} proved the  local existence and uniqueness of strong solutions to the 3-D nonhomogeneous Navier-Stokes equations for any initial data $0\leq\rho_0\in L^{\infty}(\Omega)\cap W^{1,\gamma}(\Omega),  u_0\in H^{1}_{0,\sigma}(\Omega)$, with $\gamma>1$, and if $\gamma\geq2$, then the strong solution is unique.
 More recently,  Danchin and Mucha \cite{DM}   studied  the existence and uniqueness issue for the multi-dimensional  nonhomogeneous incompressible Navier-Stokes equations supplemented with $H^{1}(\Omega)$ initial velocity and only bounded nonnegative density. Specifically, they established  the global existence of the solution for
general data in 2-D, and the global in time if the velocity satisfies a suitable
scaling-invariant smallness condition in 3-D.

Concerning the model considered in this paper, let us recall that, for the 3-D case,  under certain
assumptions,  Lukaszewicz \cite{Lu}
established the existence of weak solutions for a short time by using linearization and an almost fixed point theorem. Braz e Silva and collaborators in \cite{PB1,PB2} constructed the existence of global in time weak solutions  of  system \eqref{eq:NA} when the initial density is not necessarily
strictly positive. Local existence  of strong solutions to   system \eqref{eq:NA} was  constructed by Lukaszewicz \cite{Lu} when the initial density is strictly separated from zero.  Using a spectral semi-Galerkin method, when the initial density is bounded and away
from zero,  Boldrini et al. \cite{JB} proved the unique  local solvability  of  strong solution and some global existence results for small  data. In 2020,  Braz e Silva et al. \cite{PFM} proved the global existence and uniqueness of solution to the 3-D system \eqref{eq:NA} with  the initial velocities $(u_0,\omega_0)\in H^{1}(\mathbb{R}^{3})\times H^{1}(\mathbb{R}^{3})$  and  with  initial density being bounded from above and below by some positive constants. Obviously, the result does not allow  the presence of initial vacuum. In addition,  the  corresponding  result of the model in 2-D is also unknown.

The main  goal of this paper is to study  the global well-posedness for   $d$-dimensional ($d=2,3$) system \eqref{eq:NA} when  the initial vacuum is allowed.
 Now, let us explain some of the main difficulties and  the strategies to overcome them in the process. First, since the density is only bounded, it seems impossible to prove the uniqueness of the solution in the Eulerian coordinates as in \cite{CHK,Lijinkai}. Indeed, let $(\rho_1,u_1,\omega_1)$ and $(\rho_2,u_2,\omega_2)$ be two different solutions of system \eqref{eq:NA}. Then
 $\delta \rho=\rho_1-\rho_2$ satisfies
 $$\partial_{t}\delta \rho +u_1\cdot \nabla \delta \rho=-(u_1-u_2)\cdot\nabla \rho_2.$$
Without extra assumptions about the regularity of these solutions, the term $(u_1-u_2)\cdot\nabla \rho_2$ cannot be handled by the energy method because the usual technique to prove uniqueness via Gronwall's inequality  cannot be applied here. And the uniqueness result of Germain \cite{PG} cannot be applied here either, which requires the density function satisfying
 $\nabla \rho \in L^{\infty}\big(0,T; L^{n}(\mathbb{R}^{n})\big)$. To overcome this difficulty, the proof of uniqueness of the solution will use the Lagrangian approach, which is motivated by \cite{DM0,DH,DM}. Second, when the vacuum appears, that is, the initial density is not bounded from  below by some positive constant,  system \eqref{eq:NA} degenerates in vacuum regions  and the terms $\rho v_t, \rho \omega_t$ in the  equations are likely to vanish in some parts of the fluid domain,  which brings some difficulties for our analysis. Meanwhile,   Lemma 9 in \cite{PFM}  which plays a key role in the proof of the uniqueness  does  not apply to the presence of initial vacuum since it requires  positive lower bound of the density, which makes us have to choose another  suitable  method  different from \cite{PFM} to address the problem. Under the new framework, the estimates of  nonlinear terms
  including $\Cl u$, $\Cl\omega$ and $\nabla\Dv\omega$ in Lagrangian coordinates bring some difficulties. To overcome them, we  introduce some useful analysis tools,  for example,  the Lagrangian vorticity and the Piola identity. Third, when the density is rough and the vacuum is taken into consideration, propagate enough regularity for the velocity is the main difficulty. In most evolutionary fluid mechanics models, the uniqueness issue is
closely connected to the Lipschitz control of the flow of the velocity field. In order to achieve the
$L^{1}(0, T;W^{1,\infty})$  estimate of the velocity, the main tools are to perform the time-weighted estimates to the system \eqref{eq:NA} in the spirit of \cite{Lijinkai} and the  shifts of integrability from the time variable to the space variable in \cite{DM}. At last, much
more complicate nonlinear terms and the coupling effects of  system \eqref{eq:NA}
will also bring some troubles in our proof.

 For simplification,  we first define  the following some constants:
 \begin{align}\label{2.2-3-1-111}M\eqdefa\int_{\Omega}\rho_{0} dx, \end{align}
\begin{align}\label{2.2-3-1}C_0\eqdefa\|\sqrt{\rho_{0}}u_{0}\|_{2}^{2}+\|\sqrt{\rho_{0}}\omega_{0}\|_{2}^{2},\end{align}
\begin{align}\label{2.2-2-1A}J_0\eqdefa\nu\|\nabla u_0\|_{2}^{2}+\gamma\|\nabla \omega_0\|_{2}^{2}+\kappa\|\Dv \omega_0\|_{2}^{2}+4\chi\|\omega_0\|_{2}^{2},\end{align}
and
\begin{align}\label{2.2-2-1}K_0\eqdefa\mu\|\nabla u_0\|_{2}^{2}+\gamma\|\nabla \omega_0\|_{2}^{2}+\kappa\|\Dv \omega_0\|_{2}^{2}+\chi\|\Cl  u_0-2\omega_0\|_{2}^{2}.\end{align}
Now our main results in this paper can be  listed as follows. Let us start with
 the 2-D case.
\begin{Theorem}\label{th:main1}
Let $\Omega$ be a $C^{2}$ bounded subset of  $\mathbb{R}^{^{2}}$, or the torus $\mathbb{T}^{^{2}}$. Suppose that the initial  data $(\rho_{0},u_{0}, \omega_{0})\in L^{\infty}(\Omega)\times H_{0}^{1}(\Omega)\times H_{0}^{1}(\Omega)$ satisfy for some constant $\rho^{*}>0$,
\begin{align}\label{2.1}
0\leq \rho_{0}\leq\rho^{*},\quad \Dv u_{0}=0 \quad \text{and}\quad M>0.
\end{align}
Then system \eqref{eq:NA} has a unique globally defined solution $(\rho,u, \omega, \nabla P)$ satisfying
$$\rho\in L^{\infty}\big(\mathbb{R}^{+};L^{\infty}(\Omega)\big),\quad  u, \omega\in L^{\infty}\big(\mathbb{R}^{+};H_{0}^{1}(\Omega)\big),\quad \sqrt{\rho}u_{t},\nabla^{2}u,\nabla^{2}\omega,\nabla P\in L^{2}\big(\mathbb{R}^{+};L^{2}(\Omega)\big)$$
and also, for all $1\leq r<2 $ and $1\leq m<\infty, $
$$\nabla (\sqrt{t}P),\nabla^{2} (\sqrt{t}u),\nabla^{2} (\sqrt{t}\omega)\in L^{\infty}(0,T;L^{r}(\Omega))\cap L^{2}(0,T;L^{m}(\Omega))^{2}\quad \text{for all}\quad  T>0.$$
Furthermore, we have $\sqrt{\rho}u,\sqrt{\rho}\omega \in C\big(\mathbb{R}^{+};L^{2}(\Omega)\big), \rho\in C\big(\mathbb{R}^{+};L^{p}(\Omega)\big)$ for all  $p<\infty$, and $u,\omega\in H^{\eta}\big(0,T;L^{p}(\Omega)\big)$ for all $\eta<\frac{1}{2}$ and $T>0.$
\end{Theorem}
In the 3-D case we have the following result.
\begin{Theorem}\label{th:main2}
Let $\Omega$ be a $C^{2}$ bounded subset of $\mathbb{R}^{^{3}}$,  or the torus $\mathbb{T}^{^{3}}$. There exists a constant $\varepsilon_0>0$ such that for the initial data $(\rho_{0},u_{0},\omega_{0})\in L^{\infty}(\Omega)\times H_{0}^{1}(\Omega)\times H_{0}^{1}(\Omega)$ satisfying \eqref{2.1} and
\begin{align}\label{2.2}
 (\rho^{\ast})^{\frac{3}{2}}C_0K_0\leq \varepsilon_0.
\end{align}
Then system \eqref{eq:NA}  has a unique globally defined solution $(\rho,u,\omega,\nabla P)$ satisfying
$$\rho\in L^{\infty}\big(\mathbb{R}^{+};L^{\infty}(\Omega)\big),\quad  u,\omega\in L^{\infty}\big(\mathbb{R}^{+};H_{0}^{1}(\Omega)\big), \quad \sqrt{\rho}u_{t},\nabla^{2}u,\nabla^{2}\omega,\nabla P\in L^{2}\big(\mathbb{R}^{+};L^{2}(\Omega)\big)$$
and $$\nabla (\sqrt{t}P),\nabla^{2} (\sqrt{t}u),\nabla^{2} (\sqrt{t}\omega)\in L^{\infty}\big(0,T;L^{2}(\Omega)\big)\cap L^{2}\big(0,T;L^{6}(\Omega)\big)^{2}\quad \text{for all}\quad T>0.$$
Furthermore, we have $\sqrt{\rho}u,\sqrt{\rho}\omega\in C\big(\mathbb{R}^{+};L^{2}(\Omega)\big),\rho\in C\big(\mathbb{R}^{+};L^{p}(\Omega)\big)$ for all  $p<\infty$, and $u,\omega\in H^{\eta}\big(0,T;L^{6}(\Omega)\big)$ for all $\eta<\frac{1}{2}$ and $T>0.$
\end{Theorem}
\begin{Remark}\label{1.2} Compared with \cite{PFM},  our results allow the persistence of the vacuum and show  the global well-posedness of  solution with large data in 2-D.
\end{Remark}
\begin{Remark}\label{1.2} Here, we should point out that the fluid domain $\Omega$ is either
the   torus $\mathbb{T}^d$ or a $\cC^2$ simply connected bounded domain of $\mathbb{R}^d (d=2,3)$ since we need to use Poincar\'{e}'s inequality.
The  case $\mathbb{R}^d$ will be considered in our future work.
\end{Remark}
\noindent{\bf Notations.}  We assume $C$ be
a positive generic constant throughout this paper that may vary at
different places. By $\nabla$ we denote the gradient with respect to space variables, and by $\partial_t u$ or $u_t,$ the time derivative of function $u$.
   By $\|\cdot\|_p,$ we mean $p$-power Lebesgue norms over $\Omega$; we denote  by $H^s$ and $W^s_p$  the Sobolev
   (Slobodeckij for $s$ not an integer) space,  and put $H^s=W^s_2.$ Finally, as a great part of our analysis will concern $H^1$ regularity and will work indistinctly in
a bounded domain or in the torus, we shall adopt slightly abusively the notation $H^1_0(\Omega)$
to designate the set of $H^1(\Omega)$ functions that vanish at the boundary if $\Omega$ is a bounded domain,
or general $H^1(\mathbb{T}^d)$ functions if $\Omega=\mathbb{T}^d.$

 The rest of the paper unfolds as follows. In the next section,  we will prove the  global existence of   the solution  for system \eqref{eq:NA} and some of the time-weighted estimates on time derivatives in 2-D.
 In Section 3, we will deal with  the 3-D case.  At last, section 4 is devoted to the proof of the uniqueness of the solution to system \eqref{eq:NA}.
\par
\section{Existence of solution  and  weighted energy method in 2-D}
\subsection{Existence of solution in 2-D}
The proof is based on \emph{a priori} bounds for suitable smoothed out approximate solutions with no vacuum, then to pass to the limit. Let $j_{\delta}$ be the standard Friedrich's  mollifier and define
$$ u^{\delta}_0=j_{\delta}\ast u_0, \quad \omega^{\delta}_0=j_{\delta}\ast\omega_0,$$ and
$$ \rho^{\delta}_0=j_{\delta}\ast\rho_0,\quad \delta\leq\rho^{\delta}_0\leq\rho^{\ast}.$$
In what follows, we shall only derive  \emph{a priori} uniform energy estimates
for the approximate sequences $(\rho^{\delta},u^{\delta}, \omega^{\delta})$. Then the existence part of Theorem \ref{th:main1}
essentially follows from the \emph{a priori} estimates  and a standard compactness argument. We omit the superscript $\delta$ to keep the notation simple.

From $\rho_{t}=-u\cdot\nabla\rho$, we can easily get
\begin{equation*} \|\rho(t)\|_{L^{\infty}}=\|\rho_0\|_{L^{\infty}},\end{equation*} and with \eqref{2.1}, we have
\begin{equation}\label{3.1-1} 0\leq\rho(t)\leq\rho^{\ast},\quad (t,x)\in [0,\infty)\times \mathbb{T}^{2}.\end{equation}
Taking the $L^2$-scalar product of the
first equation of  system \eqref{eq:NA} with $u$  and
the third equation with $\omega$ respectively and the combining together,   using  $\rho_{t}=-u\cdot\nabla\rho$ and integrating it over $[0,t]$,  we obtain
\begin{equation}\label{3.1}
\begin{split}\|\sqrt{\rho}u\|_{2}^{2}+\|\sqrt{\rho}\omega\|_{2}^{2}+\int_{0}^{t}\Big(\|\nabla u\|_{2}^{2}+\|\nabla \omega\|_{2}^{2}\Big)d\tau\leq CC_0,
\end{split}
\end{equation}
where $C$ is a  universal positive constant depending the parameters $\nu, \gamma,\kappa$ and $\chi$.
\begin{Proposition}
\label{Pro:3.1}  Under the assumptions of Theorem \ref{th:main1}.
Let $(\rho,u,\omega)$ be a smooth enough solution to system \eqref{eq:NA} satisfying \eqref{3.1-1} and $T>0$. There exists a constant $C_1$ depending only on $M,\|\rho_{0}\|_{2},\|\sqrt{\rho_{0}}u_{0}\|_{2}$ and  $\rho^{*}$ so that for all $t\in[0,T)$, we have
\begin{equation}\label{3.2-1}
\begin{split}
&\|\nabla u\|_{2}^{2}+\|\nabla \omega\|_{2}^{2}+\int_{0}^{t}\Big(\|\sqrt{\rho}u_{t}\|_{2}^{2}+\|\sqrt{\rho}\omega_{t}\|_{2}^{2}+\|\nabla^{2}u\|_{2}^{2}+\|\nabla^{2}\omega\|_{2}^{2}+\|\nabla P\|^{2}_{2}\Big)d\tau\\
&\leq \Big(e+CJ_0\Big)^{\exp\big(CC_1C_{0}\big)}.
\end{split}
\end{equation}
Furthermore, for all $p\in[1,\infty)$ and $t\in[0,T)$, we have
\begin{align}\label{3.3A}
\|u\|_{p}+\|\omega\|_{p}\leq C\frac{C_{0}}{M}+C_{p}\Big(1+\frac{\|M-\rho\|_{2}}{M}\Big)\Big(\|\nabla u\|_{2}+\|\nabla \omega\|_{2}\Big).
\end{align}
\end{Proposition}
\noindent{\bf Proof.}\   Taking the $L^2$-scalar product of the
first equation of  system \eqref{eq:NA} with $u_{t}$  and
the third equation with $\omega_{t}$ respectively, we obtain from  $\rho_{t}=-u\cdot\nabla\rho$ that 
\begin{equation*}
\int_{\mathbb{T}^{2}}\rho| u_{t}|^{2}dx+\frac{\nu}{2}\frac{d}{dt}\int_{\mathbb{T}^{2}}|\nabla u|^{2}dx
=2\chi\int_{\mathbb{T}^{2}}\Cl \omega\cdot u_{t}dx-\int_{\mathbb{T}^{2}}(\rho u\cdot\nabla u)\cdot u_{t}dx,
\end{equation*}
and
\begin{equation*}
\begin{split}
\int_{\mathbb{T}^{2}}\rho| \omega_{t}|^{2}dx+&\frac{\gamma}{2}\frac{d}{dt}\int_{\mathbb{T}^{2}}|\nabla \omega|^{2}dx
+\frac{\kappa}{2}\frac{d}{dt}\int_{\mathbb{T}^{2}}|\Dv \omega|^{2}dx
+2\chi\frac{d}{dt}\int_{\mathbb{T}^{2}}|\omega|^{2}dx\\
&=2\chi\int_{\mathbb{T}^{2}}\Cl u\cdot \omega_{t}dx-\int_{\mathbb{T}^{2}}(\rho u\cdot\nabla \omega)\cdot \omega_{t}dx.
\end{split}
\end{equation*}
Adding the two identities and using H\"{o}lder's and Young's inequalities yield that
\begin{equation*}
\begin{split}
&\|\sqrt{\rho}u_{t}\|_{2}^{2}+\|\sqrt{\rho}\omega_{t}\|_{2}^{2}
+\frac{1}{2}\frac{d}{dt}\Big(\nu\|\nabla u\|_{2}^{2}+\gamma\|\nabla \omega\|_{2}^{2}+\kappa\|\Dv \omega\|_{2}^{2}+4\chi\|\omega\|_{2}^{2}+4\chi(\Cl  u,\omega)\Big)\\
&=-\int_{\mathbb{T}^{2}}(\rho u\cdot\nabla u)\cdot u_{t}dx-\int_{\mathbb{T}^{2}}(\rho u\cdot\nabla \omega)\cdot \omega_{t}dx\\
&\leq\frac{1}{2}\int_{\mathbb{T}^{2}}\rho|u_{t}|^{2}dx+\frac{1}{2}\int_{\mathbb{T}^{2}}\rho|\omega_{t}|^{2}dx
+\frac{1}{2}\int_{\mathbb{T}^{2}}\rho|u\cdot\nabla u|^{2}dx
+\frac{1}{2}\int_{\mathbb{T}^{2}}\rho|u\cdot\nabla \omega|^{2}dx,
\end{split}
\end{equation*}
which implies that
\begin{equation}\label{3.5-C}
\begin{split}
&\|\sqrt{\rho}u_{t}\|_{2}^{2}+\|\sqrt{\rho}\omega_{t}\|_{2}^{2}
+\frac{d}{dt}\Big(\nu\|\nabla u\|_{2}^{2}+\gamma\|\nabla \omega\|_{2}^{2}+\kappa\|\Dv \omega\|_{2}^{2}+4\chi\|\omega\|_{2}^{2}+4\chi(\Cl  u,\omega)\Big)\\
&\leq\int_{\mathbb{T}^{2}}\rho|u\cdot\nabla u|^{2}dx
+\int_{\mathbb{T}^{2}}\rho|u\cdot\nabla \omega|^{2}dx.
\end{split}
\end{equation}
Due to \begin{equation*}
\begin{split}
4\chi\big|(\Cl  u,\omega)\big|\leq \chi\|\nabla u\|^2_{2}+4\chi\|\omega\|^2_{2},
\end{split}
\end{equation*}
then there exist two positive constants $c_1$ and $c_2$ such that
\begin{equation*}
\begin{split}
c_1\alpha_1\leq\nu\|\nabla u\|_{2}^{2}+\gamma\|\nabla \omega\|_{2}^{2}+\kappa\|\Dv \omega\|_{2}^{2}+4\chi\|\omega\|_{2}^{2}+4\chi(\Cl  u,\omega)
\leq c_2\alpha_1
\end{split}
\end{equation*}
with $\alpha_1\eqdefa\nu\|\nabla u\|_{2}^{2}+\gamma\|\nabla \omega\|_{2}^{2}+\kappa\|\Dv \omega\|_{2}^{2}+4\chi\|\omega\|_{2}^{2}.$\\
This along with \eqref{3.5-C} ensures that 
\begin{equation}\label{3.5-B}
\begin{split}
&\|\sqrt{\rho}u_{t}\|_{2}^{2}+\|\sqrt{\rho}\omega_{t}\|_{2}^{2}
+\frac{d}{dt}\Big(\nu\|\nabla u\|_{2}^{2}+\gamma\|\nabla \omega\|_{2}^{2}+\kappa\|\Dv \omega\|_{2}^{2}+4\chi\|\omega\|_{2}^{2}\Big)\\
&\leq C\Big(\int_{\mathbb{T}^{2}}\rho|u\cdot\nabla u|^{2}dx
+\int_{\mathbb{T}^{2}}\rho|u\cdot\nabla \omega|^{2}dx\Big).
\end{split}
\end{equation}
It follows from integrating with respect to time over $[0,t]$   on the both sides of \eqref{3.5-B} that
\begin{equation}\label{3.5-1A}
\begin{split}
&\|\nabla u\|_{2}^{2}+\|\nabla \omega\|_{2}^{2}+\int_{0}^{t}\Big(\|\sqrt{\rho}u_{t}\|_{2}^{2}+\|\sqrt{\rho}\omega_{t}\|_{2}^{2}\Big)d\tau\\
&\leq CJ_0+C\Big(\int_{0}^{t}\int_{\mathbb{T}^{2}}\rho|u\cdot\nabla u|^{2}dxd\tau
+\int_{0}^{t}\int_{\mathbb{T}^{2}}\rho|u\cdot\nabla \omega|^{2}dxd\tau\Big),
\end{split}
\end{equation}
where  $J_0$ is given by \eqref{2.2-2-1A}.

In order to bound the second derivatives of $(u,\omega)$, let us take the  $L^2$-scalar product of the
first equation of  system \eqref{eq:NA} with $-\Delta u$  and
the third equation with $-\Delta\omega$ respectively.  Adding the resulting equations, we easily arrive at
\begin{equation*}
\begin{split}
&\nu\|\nabla^{2}u\|_{2}^{2}+\gamma\|\nabla^{2}\omega\|_{2}^{2}+\kappa\|\nabla\Dv \omega\|_{2}^{2}+4\chi\|\nabla\omega\|_{2}^{2}\\
&=4\chi(\Cl \omega,-\Delta u)+(\rho u_t,\Delta u)+(\rho \omega_t,\Delta \omega)\\
&\quad+(\rho u\cdot\nabla u,\Delta u)+(\rho u\cdot\nabla \omega,\Delta \omega),
\end{split}
\end{equation*}
where we have used the fact $(\Cl \omega,-\Delta \omega)=(\Cl \omega,-\Delta u)$.\\
Hence, using  Young's inequality yields that
\begin{equation}\label{3.6-1}
\begin{split}
\|\nabla^{2}u\|_{2}^{2}+\|\nabla^{2}\omega\|_{2}^{2}\leq C\Big(\|\sqrt{\rho}u_{t}\|_{2}^{2}+\|\sqrt{\rho}\omega_{t}\|_{2}^{2}
+\int_{\mathbb{T}^{2}}\rho|u\cdot\nabla u|^{2}dx
+\int_{\mathbb{T}^{2}}\rho|u\cdot\nabla \omega|^{2}dx \Big).
\end{split}
\end{equation}
Next, we deal with  the gradient of the pressure. Taking the divergence of the linear momentum equation gives
\begin{equation*}\Delta P=-\D\big(\rho u_t+\rho u\cdot\nabla u\big),
\end{equation*} and consequently the pressure $P$ may be recovered by
\begin{equation*} \nabla P=-\nabla\Delta^{-1}\Dv\big(\rho u_t+\rho u\cdot\nabla u\big).
\end{equation*}
By the bounded of Riesz's  operator, we have
\begin{equation}\label{3.6-2} \|\nabla P\|^{2}_{2}\leq C\Big(\|\sqrt{\rho}u_{t}\|_{2}^{2}+\int_{\mathbb{T}^{2}}\rho|u\cdot\nabla u|^{2}dx
\Big).
\end{equation}
We finally conclude from \eqref{3.5-1A}, \eqref{3.6-1} and \eqref{3.6-2},     that
\begin{equation}\label{3.6-3A}
\begin{split}
&\|\nabla u\|_{2}^{2}+\|\nabla \omega\|_{2}^{2}+\int_{0}^{t}\Big(\|\sqrt{\rho}u_{t}\|_{2}^{2}+\|\sqrt{\rho}\omega_{t}\|_{2}^{2}+\|\nabla^{2}u\|_{2}^{2}+\|\nabla^{2}\omega\|_{2}^{2}+\|\nabla P\|^{2}_{2}\Big)d\tau\\
&\leq CJ_0+C\Big(\int_{0}^{t}\int_{\mathbb{T}^{2}}\rho|u\cdot\nabla u|^{2}dxd\tau
+\int_{0}^{t}\int_{\mathbb{T}^{2}}\rho|u\cdot\nabla \omega|^{2}dxd\tau\Big).
\end{split}
\end{equation}
Here and in what follows, we  bound the last   term  on the right-hand side of the inequality above. Using H\"{o}lder's inequality and  2-D Ladyzhenskaya's inequality \cite{OL}:
$\|v\|^{2}_{4}\leq C\|v\|_{2}\|\nabla v\|_{2}$, we deduce that
\begin{equation}\label{3.8-1}
\int_{\mathbb{T}^{2}}\rho|u\cdot\nabla u|^{2}dx\leq\|\sqrt{\rho}|u|^{2}\|_{2}\|\sqrt{\rho}|\nabla u|^{2}\|_{2}\leq C\sqrt{\rho^{*}}\|\sqrt{\rho}|u|^{2}\|_{2}\|\nabla u\|_{2}\|\nabla^{2} u\|_{2},
\end{equation} and
\begin{equation}\label{3.8-2}
\int_{\mathbb{T}^{2}}\rho|u\cdot\nabla \omega|^{2}dx\leq\|\sqrt{\rho}|u|^{2}\|_{2}\|\sqrt{\rho}|\nabla \omega|^{2}\|_{2}\leq C\sqrt{\rho^{*}}\|\sqrt{\rho}|u|^{2}\|_{2}\|\nabla \omega\|_{2}\|\nabla^{2} \omega\|_{2}.
\end{equation}
For the term $\|\sqrt{\rho}|u|^{2}\|_{2}^{2}$, we have
\begin{equation*}\label{}
\|\sqrt{\rho}|u|^{2}\|_{2}^{2}\leq C\|\sqrt{\rho}u\|_{2}^{2}\|\nabla u\|_{2}^{2}\log\Big(e+\frac{\|\rho_{0}-M\|_{2}^{2}}{M^{2}}+\frac{\rho^{*}\|\nabla u\|_{2}^{2}}{\|\sqrt{\rho}u\|_{2}^{2}}\Big),
\end{equation*}
where we have used the following improvement of  Ladyzhenskaya's
inequality that has been pointed out by B. Desjardins in \cite{BD}:
\begin{equation*}\label{}
\|\sqrt{\rho}v^{2}\|_{2}\leq C\|\sqrt{\rho}v\|_{2}\|\nabla v\|_{2}\log^{\frac{1}{2}}\Big(e+\frac{\|\rho-M\|_{2}^{2}}{M^{2}}+\frac{\rho^{*}\|\nabla v\|_{2}^{2}}{\|\sqrt{\rho}v\|_{2}^{2}}\Big),
\end{equation*}
where for all $v\in H^{1}(\mathbb{T}^{2})$ and $\rho\in L^{\infty}(\mathbb{T}^{2})$ with $0\leq\rho\leq\rho^{*}$. \\
Observing that the function $z\to z\lg(e+\frac{1}{z})$ is increasing,  and reverting to \eqref{3.8-1} and \eqref{3.8-2}, we end up with
\begin{equation}\label{3.16-1-AA}
\begin{split}
\int_{\mathbb{T}^{2}}\rho|u\cdot\nabla u|^{2}dx
&\leq\varepsilon\|\nabla^{2} u\|_{2}^{2}+C_{\rho^{*}}\|\sqrt{\rho} |u|^{2}\|_{2}^{2}\|\nabla u\|_{2}^{2}\\
&\leq\varepsilon\|\nabla^{2} u\|_{2}^{2}+C_{\rho^{*}}C_{0}\|\nabla u\|_{2}^{4}\lg\Big(e+\frac{\|\rho_{0}-M\|_{2}^{2}}{M^{2}}+\frac{\rho^{*}\|\nabla u\|_{2}^{2}}{C_0}\Big),
\end{split}
\end{equation}
and
\begin{equation}\label{3.16-1-BB}
\begin{split}
\int_{\mathbb{T}^{2}}\rho|u\cdot\nabla \omega|^{2}dx
&\leq\varepsilon\|\nabla^{2} \omega\|_{2}^{2}+C_{\rho^{*}}\|\sqrt{\rho} |u|^{2}\|_{2}^{2}\|\nabla \omega\|_{2}^{2}\\
&\leq\varepsilon\|\nabla^{2} \omega\|_{2}^{2}+C_{\rho^{*}}C_{0}\|\nabla u\|_{2}^{2}\|\nabla \omega\|_{2}^{2}\lg\Big(e+\frac{\|\rho_{0}-M\|_{2}^{2}}{M^{2}}+\frac{\rho^{*}\|\nabla u\|_{2}^{2}}{C_{0}}\Big),
\end{split}
\end{equation}
where $\varepsilon$ is arbitrary small positive constant.\\
Then combining  \eqref{3.16-1-AA} and \eqref{3.16-1-BB} with \eqref{3.6-3A} yields that
\begin{equation}\label{3.16-1}
\begin{split}
&\|\nabla u\|_{2}^{2}+\|\nabla \omega\|_{2}^{2}+\int_{0}^{t}\Big(\|\sqrt{\rho}u_{t}\|_{2}^{2}+\|\sqrt{\rho}\omega_{t}\|_{2}^{2}+\|\nabla^{2}u\|_{2}^{2}+\|\nabla^{2}\omega\|_{2}^{2}+\|\nabla P\|^{2}_{2}\Big)d\tau\\
&\leq CJ_0+CC_1\Big(\int_{0}^{t}\lg(e+\|\nabla u\|_{2}^{2})\|\nabla u\|_{2}^{2}(\|\nabla u\|_{2}^{2}+\|\nabla \omega\|_{2}^{2})d\tau\Big)
\end{split}
\end{equation}
with $C_{1}$ depending only on $\rho^{*}, C_{0}, M$ and $\|\rho_{0}\|_{2}.$\\
Denoting $f(t):=C_{1}\|\nabla u\|_{2}^{2}$ and
\begin{equation*}X(t):=\|\nabla u\|_{2}^{2}+\|\nabla \omega\|_{2}^{2}+\int_{0}^{t}\Big(\|\sqrt{\rho}u_{t}\|_{2}^{2}+\|\sqrt{\rho}\omega_{t}\|_{2}^{2}+\|\nabla^{2}u\|_{2}^{2}+\|\nabla^{2}\omega\|_{2}^{2}+\|\nabla P\|^{2}_{2}\Big)d\tau.\end{equation*}
It follows from the inequality \eqref{3.16-1}   that
\begin{equation}\label{3.16-2}X(t)\leq CJ_0 +C\int_{0}^{t} f(\tau)X(\tau)\lg\Big(e+X(\tau)\Big)d\tau.\end{equation}
Setting $g(t)=\int_{0}^{t} f(\tau)X(\tau)\lg\big(e+X(\tau)\big)d\tau$,  from \eqref{3.16-2}, we obtain
$e+X(t)\leq  e+CJ_0 +Cg(t).$  Then,
\begin{equation*}
\begin{split}\frac{d}{dt}g(t)&= f(t)X(t)\lg\Big(e+X(t)\Big)\\
&\leq f(t) \Big(e+X(t)\Big)\lg\big(e+X(t)\Big)\\
&\leq f(t) \Big(e+CJ_0 +Cg(t)\Big)\lg\Big( e+CJ_0 +Cg(t)\Big),
\end{split}
\end{equation*}
from which we get, for all $t\geq0,$
$$e+CJ_0 +Cg(t)\leq\Big(e+CJ_0\Big)^{\exp\Big(C\int_{0}^{t}f(\tau)d\tau\Big)}.$$
Hence,
\begin{equation*}
\begin{split}
&\|\nabla u\|_{2}^{2}+\|\nabla \omega\|_{2}^{2}+\int_{0}^{t}\Big(\|\sqrt{\rho}u_{t}\|_{2}^{2}+\|\sqrt{\rho}\omega_{t}\|_{2}^{2}+\|\nabla^{2}u\|_{2}^{2}+\|\nabla^{2}\omega\|_{2}^{2}+\|\nabla P\|^{2}_{2}\Big)d\tau\\
&\leq \Big(e+CJ_0\Big)^{\exp\big(CC_1C_{0}\big)}.
\end{split}
\end{equation*}

In order to prove \eqref{3.3A},  denoting by $\bar{u}(t), \bar{\omega}(t)$ the average of $u(t), \omega(t)$ on $\mathbb{T}^{2}$,  for all $p\in[1,\infty),$  it then  follows from  Sobolev embedding that
\begin{equation}\label{3.10}
\begin{split}
\|u(t)\|_{p}+\|\omega(t)\|_{p}
&\leq|\bar{u}(t)|+|\bar{w}(t)|+\|u(t)-\bar{u}(t)\|_{p}+\|\omega(t)-\bar{\omega}(t)\|_{p}\\
&\leq|\bar{u}(t)|+|\bar{w}(t)|+C_{p}\Big(\|\nabla u(t)\|_{2}+\|\nabla\omega(t)\|_{2}\Big).
\end{split}
\end{equation}
On the other hand, applying Poincar\'{e}'s inequality yields that
\begin{equation*}
\begin{split}
M\big(|\bar{u}(t)|+|\bar{\omega}(t)|\big)&=\int_{\mathbb{T}^{2}}\rho udx+\int_{\mathbb{T}^{2}}\rho \omega dx+
\int_{\mathbb{T}^{2}}\big(M-\rho\big)\big(u-\bar{u}\big)dx\\&\quad+\int_{\mathbb{T}^{2}}\big(M-\rho\big)\big(\omega-\bar{\omega}\big)dx\\
&\leq CC_0+\|M-\rho\|_{2}\Big(\|\nabla u\|_{2}+\|\nabla \omega\|_{2}\Big).
\end{split}
\end{equation*}
Putting that latter inequality into \eqref{3.10} yields \eqref{3.3A}.
\subsection{Weighted energy method in 2-D.}
In order to obtain the shift integrability from time to
space variables,  our next aim is to exploit some bounds, for example, $(\sqrt{\rho t}u_{t},\sqrt{\rho t}\omega_{t})$ in $L^{\infty}([0,T];L^{2})$  and $(\sqrt{t}\nabla u_{t},\sqrt{t}\nabla \omega_{t})$ in $L^{2}([0,T];L^{2})$ respectively, in terms of the data.
\begin{Lemma}\label{3.3-1}
Assume $d=2$ and that the solution is smooth enough of  system \eqref{eq:NA} with no vacuum. Then for all $t\geq0,$ we have
\begin{align}\label{3.21}
\|\sqrt{\rho t}u_{t}\|_{2}^{2}+\|\sqrt{\rho t}\omega_{t}\|_{2}^{2}
+\int_{0}^{t}\tau\|\nabla u_{t}\|_{2}^{2}d\tau d\tau
+\int_{0}^{t}\tau\|\nabla\omega_{t}\|_{2}^{2}
\leq \exp\Big(\int_{0}^{t}h_1(\tau)d\tau\Big)-1
\end{align}
with $h_1\in L^{1}_{loc}(\mathbb{R}^{+})$ depending only on $\rho^{*},\|\sqrt{\rho_{0}}u_{0}\|_{2},
\|\sqrt{\rho_{0}}\omega_{0}\|_{2}$ and $K_{0}$.
\end{Lemma}
\noindent{\bf Proof.}\ Differentiating $\eqref{eq:NA}_{1}$ and $\eqref{eq:NA}_{3}$ with respect to $t$, respectively, we have
$$\rho u_{tt}+\rho_{t}u_{t}+\rho_{t}u\cdot\nabla u+\rho u_{t}\cdot\nabla u+\rho u\cdot\nabla u_{t}-\nu\Delta u_{t}+\nabla P_{t}=2\chi \Cl \omega_{t},$$
and
$$\rho \omega_{tt}+\rho_{t}\omega_{t}+\rho_{t}u\cdot\nabla \omega+\rho u_{t}\cdot\nabla \omega+\rho u\cdot\nabla \omega_{t}-\gamma\Delta \omega_{t}-\kappa\nabla\Dv\omega_{t}+4\chi\omega_{t}
=2\chi \Cl u_{t}.$$
Then, multiplying by $\sqrt{t}$  the above two inequalities respectively yields that 
\begin{equation}\label{2.19-A}
\begin{split}\rho (\sqrt{t}u_{t})_{t}-&\frac{1}{2\sqrt{t}}\rho u_{t}+\sqrt{t}\rho_{t}u_{t}
+\sqrt{t}\rho_{t}u\cdot\nabla u+\sqrt{t}\rho u_{t}\cdot\nabla u+\sqrt{t}\rho u\cdot\nabla u_{t}-\nu\Delta (\sqrt{t}u_{t})+\nabla (\sqrt{t}P_{t})\\&=2\chi \Cl (\sqrt{t}\omega_{t}),\end{split}
\end{equation}
and
\begin{equation}\label{2.20-A}
\begin{split}
\rho (\sqrt{t}\omega_{t})_{t}&-\frac{1}{2\sqrt{t}}\rho \omega_{t}
+\sqrt{t}\rho_{t}\omega_{t}+\sqrt{t}\rho_{t}u\cdot\nabla \omega
+\sqrt{t}\rho u_{t}\cdot\nabla \omega+\sqrt{t}\rho u\cdot\nabla \omega_{t}-\gamma\Delta (\sqrt{t}\omega_{t})\\
&-\kappa\nabla\Dv(\sqrt{t}\omega_{t})+4\chi(\sqrt{t}\omega_{t})
=2\chi \Cl (\sqrt{t}u_{t}).
\end{split}
\end{equation}
Taking the $L^{2}$ scalar product of \eqref{2.19-A} with $\sqrt{t}u_{t}$ and \eqref{2.20-A} with $\sqrt{t}\omega_{t}$ respectively, we get
\begin{equation}\label{3.11-AAA}
\begin{split}
&\frac{1}{2}\frac{d}{dt}\int_{\mathbb{T}^{2}}\rho t|u_{t}|^{2}dx+\nu\int_{\mathbb{T}^{d}}t|\nabla u_{t}|^{2}dx\\
&\leq\frac{1}{2}\int_{\mathbb{T}^{2}}\rho|u_{t}|^{2}dx-\frac{1}{2}\int_{\mathbb{T}^{d}}t\rho_{t}
|u_{t}|^{2}dx\\
&\quad-\int_{\mathbb{T}^{2}}(\sqrt{t}\rho_{t}u\cdot\nabla u)\cdot(\sqrt{t}u_{t})dx-
\int_{\mathbb{T}^{2}}(\sqrt{t}\rho u_{t}\cdot\nabla u)\cdot(\sqrt{t}u_{t})dx
\\
&\quad-\int_{\mathbb{T}^{2}}(\sqrt{t}\rho u\cdot\nabla u_{t})\cdot(\sqrt{t}u_{t})dx
+2\chi\int_{\mathbb{T}^{2}}\Cl(\sqrt{t}\omega_{t})\cdot(\sqrt{t}u_{t})dx,
\end{split}
\end{equation}
and
\begin{equation}\label{3.11-AAAA}
\begin{split}
&\frac{1}{2}\frac{d}{dt}\int_{\mathbb{T}^{2}}\rho t|\omega_{t}|^{2}dx
+\gamma\int_{\mathbb{T}^{2}}t|\nabla \omega_{t}|^{2}dx
+4\chi\int_{\mathbb{T}^{2}}t|\omega_{t}|^{2}dx\\
&\leq\frac{1}{2}\int_{\mathbb{T}^{2}}\rho|\omega_{t}|^{2}dx-\frac{1}{2}\int_{\mathbb{T}^{d}}t\rho_{t}
|\omega_{t}|^{2}dx\\
&\quad-\int_{\mathbb{T}^{2}}(\sqrt{t}\rho_{t}u\cdot\nabla \omega)\cdot(\sqrt{t}\omega_{t})dx-\int_{\mathbb{T}^{2}}(\sqrt{t}\rho u_{t}\cdot\nabla \omega)\cdot(\sqrt{t}\omega_{t})dx\\
&\quad-
\int_{\mathbb{T}^{2}}(\sqrt{t}\rho u\cdot\nabla \omega_{t})\cdot(\sqrt{t}\omega_{t})dx
+2\chi\int_{\mathbb{T}^{2}}\Cl(\sqrt{t}u_{t})\cdot(\sqrt{t}\omega_{t})dx.
\end{split}
\end{equation}
Note that
\begin{equation*}
\begin{split}
&2\chi\int_{\mathbb{T}^{2}}\Cl(\sqrt{t}\omega_{t})\cdot(\sqrt{t}u_{t})dx+
2\chi\int_{\mathbb{T}^{2}}\Cl(\sqrt{t}u_{t})\cdot(\sqrt{t}\omega_{t})dx\\
&=4\chi\int_{\mathbb{T}^{2}}\Cl(\sqrt{t}u_{t})\cdot(\sqrt{t}\omega_{t})dx\\
&\leq4\chi\|\sqrt{t}\nabla u_{t}\|_2\|\sqrt{t}\omega_{t}\|_2\\
&\leq\chi\|\sqrt{t}\nabla u_{t}\|_2^2+4\chi\|\sqrt{t}\omega_{t}\|_2^2,
\end{split}
\end{equation*}
which together with \eqref{3.11-AAA} and \eqref{3.11-AAAA}  implies that
\begin{align}\label{3.11}\begin{split}
\frac{d}{dt}&\Big(\|\sqrt{\rho t}u_{t}\|_{2}^{2}+\|\sqrt{\rho t}\omega_{t}\|_{2}^{2}\Big)dx+\int_{\mathbb{T}^{2}}t|\nabla u_{t}|^{2}dx+\int_{\mathbb{T}^{2}}t|\nabla \omega_{t}|^{2}dx\\
&\leq \int_{\mathbb{T}^{2}}\big(\rho|u_{t}|^{2}+\rho|\omega_{t}|^{2}\big)dx+\int_{\mathbb{T}^{2}}\big(t\rho_{t}
|u_{t}|^{2}+t\rho_{t}|\omega_{t}|^{2}\big)dx\\
&\quad-\int_{\mathbb{T}^{2}}\big((\sqrt{t}\rho u_{t}\cdot\nabla u)\cdot(\sqrt{t}u_{t})+(\sqrt{t}\rho u_{t}\cdot\nabla \omega)\cdot(\sqrt{t}\omega_{t})\big)dx\\
&\quad-\int_{\mathbb{T}^{2}}\big((\sqrt{t}\rho u_{t}\cdot\nabla u)\cdot(\sqrt{t}u_{t})+(\sqrt{t}\rho u_{t}\cdot\nabla \omega)\cdot(\sqrt{t}\omega_{t})\big)dx\\
&\quad-\int_{\mathbb{T}^{2}}\big((\sqrt{t}\rho u\cdot\nabla u_{t})\cdot(\sqrt{t}u_{t})
+(\sqrt{t}\rho u\cdot\nabla \omega_{t})\cdot(\sqrt{t}\omega_{t})\big)dx
\\&\triangleq \sum_{i=1}^{5}I_{i}.
\end{split}
\end{align}
In what follows, we estimate term by term above.
 For $I_{2}$,  thanks to $\rho_{t}=-u\cdot\nabla\rho$, we have
\begin{equation}\label{2.27-A-A}
\begin{split}
I_{2}&\leq C\Big|\int_{\mathbb{T}^{2}}t\Dv(\rho u)|u_{t}|^{2}dx+\int_{\mathbb{T}^{2}}t\Dv(\rho u)|\omega_{t}|^{2}dx\Big|\\
&\leq C\int_{\mathbb{T}^{2}}t\rho|u||\nabla u_{t}||u_{t}|dx+
\int_{\mathbb{T}^{2}}t\rho|u||\nabla \omega_{t}||\omega_{t}|dx\\
&\leq C\Big(\int_{\mathbb{T}^{2}}\rho t|u_{t}|^{2}dx\Big)^{\frac{1}{2}}
\Big(\int_{\mathbb{T}^{2}}t\rho|u|^{2}|\nabla u_{t}|^{2}dx\Big)^{\frac{1}{2}}
\\&\quad+C\Big(\int_{\mathbb{T}^{2}}\rho t|\omega_{t}|^{2}dx\Big)^{\frac{1}{2}}
\Big(\int_{\mathbb{T}^{2}}t\rho |u|^{2}|\nabla\omega_{t}|^{2}dx\Big)^{\frac{1}{2}}\\
&\leq C\|\sqrt{\rho t}u_{t}\|_{2}\|u\|_{\infty}\|\sqrt{t}\nabla u_{t}\|_{2}+
C\|\sqrt{\rho t}\omega_{t}\|_{2}\|u\|_{\infty}\|\sqrt{t}\nabla \omega_{t}\|_{2}\\
&\leq\varepsilon
\Big(\|\sqrt{t}\nabla u_{t}\|^{2}_{2}+\|\sqrt{t}\nabla \omega_{t}\|^{2}_{2}\Big)
+C\|u\|^{2}_{\infty}\Big(\|\sqrt{\rho t}u_{t}\|^{2}_{2}+
\|\sqrt{\rho t}\omega_{t}\|^{2}_{2}\Big).
\end{split}
\end{equation}
For $I_{3}$, according to $\rho_{t}=-u\cdot\nabla\rho$ and then  performing an integration by parts, we get
\begin{equation}\label{2.27-A}
\begin{split}
I_{3}&\leq\Big|-\int_{\mathbb{T}^{2}}(\sqrt{t}\rho_{t}u\cdot\nabla u)\cdot(\sqrt{t}u_{t})dx-
\int_{\mathbb{T}^{2}}(\sqrt{t}\rho_{t}u\cdot\nabla \omega)\cdot(\sqrt{t}\omega_{t})dx\Big|\\
&\leq\Big|-\int_{\mathbb{T}^{2}}t\rho u\cdot\nabla[(u\cdot\nabla u)\cdot u_{t}]dx
-\int_{\mathbb{T}^{2}}t\rho u\cdot\nabla[(u\cdot\nabla \omega)\cdot \omega_{t}]dx\Big|\\
&\leq\int_{\mathbb{T}^{2}}t\rho|u|\Big(|\nabla u|^{2}|u_{t}|+|u||\nabla^{2} u||u_{t}|
+|u||\nabla u||\nabla u_{t}|\\
&\quad+|\nabla u||\nabla \omega||\omega_{t}|+|u||\nabla^{2} \omega||\omega_{t}|
+|u||\nabla \omega||\nabla \omega_{t}|\Big)dx\\
&\triangleq\sum_{i=1}^{6}I_{3i}.
\end{split}
\end{equation}
It follows from   H\"{o}lder's  and Young's inequalities that 
\begin{equation*}
I_{31}=\int_{\mathbb{T}^{2}}\sqrt{\rho t}|u||\nabla u||\nabla u||\sqrt{\rho t}u_{t}|dx
\leq\|u\|_{\infty}^{2}\|\sqrt{\rho t}u_{t}\|_{2}^{2}+CT\rho^{*}\|\nabla u\|_{4}^{4},
\end{equation*}
\begin{equation*}
I_{34}=\int_{\mathbb{T}^{2}}\sqrt{\rho t}|u||\nabla u||\nabla \omega||\sqrt{\rho t}\omega_{t}|dx
\leq\|u\|_{\infty}^{2}\|\sqrt{\rho t}\omega_{t}\|_{2}^{2}+CT\rho^{*}\Big(\|\nabla u\|_{4}^{4}+\|\nabla \omega\|_{4}^{4}\Big).
\end{equation*}
Along the same line, we have
$$I_{32}=\int_{\mathbb{T}^{2}}t\rho|u|^{2}|\nabla^{2} u||u_{t}|dx
\leq \rho^{*}T\|\nabla^{2}u\|_{2}^{2}+\|u\|_{\infty}^{4}\|\sqrt{\rho t}u_{t}\|_{2}^{2},$$
$$I_{35}=\int_{\mathbb{T}^{2}}t\rho|u|^{2}|\nabla^{2} \omega||\omega_{t}|dx
\leq \rho^{*}T\|\nabla^{2}\omega\|_{2}^{2}+\|u\|_{\infty}^{4}\|\sqrt{\rho t}\omega_{t}\|_{2}^{2},$$
\begin{equation*}
\begin{split}
I_{33}=\int_{\mathbb{T}^{2}}t\rho|u|^{2}|\nabla u||\nabla u_{t}|dx
&\leq \varepsilon\int_{\mathbb{T}^{2}}|\nabla\sqrt{t}u_{t}|^{2}dx +C\int_{\mathbb{T}^{2}}t\rho^{2}|u|^{4}|\nabla u|^{2}dx\\
&\leq \varepsilon\int_{\mathbb{T}^{2}}|\nabla\sqrt{t}u_{t}|^{2}dx+
C_{T,\rho^{*}}\|u\|_{\infty}^{4}\|\nabla u\|_{2}^{2},
\end{split}
\end{equation*}
\begin{equation*}
\begin{split}
I_{36}=\int_{\mathbb{T}^{2}}t\rho|u|^{2}|\nabla \omega||\nabla\omega_{t}|dx
&\leq \varepsilon\int_{\mathbb{T}^{2}}|\nabla\sqrt{t}\omega_{t}|^{2}dx
+C\int_{\mathbb{T}^{2}}t\rho^{2}|u|^{4}|\nabla \omega|^{2}dx\\
&\leq \varepsilon\int_{\mathbb{T}^{2}}|\nabla\sqrt{t}\omega_{t}|^{2}dx+
C_{T,\rho^{*}}\|u\|_{\infty}^{4}\|\nabla\omega\|_{2}^{2}.
\end{split}
\end{equation*}
For $I_{4}$ and $I_{5}$, using  H\"{o}lder's  and Young's inequalities,  we deduce that
\begin{equation}
\begin{split}
I_{4}&\leq\|\nabla u\|_{2}\|\sqrt{\rho t}u_{t}\|_{4}^{2}+
\|\nabla\omega\|_{2}\|\sqrt{\rho t}u_{t}\|_{4}\|\sqrt{\rho t}\omega_{t}\|_{4}\\
&\leq(\rho^{*})^{3/4}\|\nabla u\|_{2}\|\sqrt{\rho t}u_{t}\|_{2}^{\frac{1}{2}}\|\sqrt{ t}u_{t}\|_{6}^{\frac{3}{2}}+
(\rho^{*})^{3/4}\|\nabla\omega\|_{2}\|\sqrt{\rho t}u_{t}\|_{2}^{\frac{1}{4}}
\|\sqrt{t}u_{t}\|_{6}^{\frac{3}{4}}\|\sqrt{\rho t}\omega_{t}\|_{2}^{\frac{1}{4}}
\|\sqrt{t}\omega_{t}\|_{6}^{\frac{3}{4}}\\
&\leq C(\rho^{*})^{3/4}\|\nabla u\|_{2}\|\sqrt{\rho t}u_{t}\|_{2}^{\frac{1}{2}}\|\sqrt{ t}\nabla u_{t}\|_{2}^{\frac{3}{2}}+
C(\rho^{*})^{3/4}\|\nabla\omega\|_{2}\|\sqrt{\rho t}u_{t}\|_{2}^{\frac{1}{4}}
\|\nabla\sqrt{t}u_{t}\|_{2}^{\frac{3}{4}}\|\sqrt{\rho t}\omega_{t}\|_{2}^{\frac{1}{4}}
\|\nabla\sqrt{t}\omega_{t}\|_{2}^{\frac{3}{4}}\\
&\leq\varepsilon\Big(\|\sqrt{t}\nabla u_{t}\|_{2}^{2}+\|\sqrt{t}\nabla\omega_{t}\|_{2}^{2}\Big)
+C_{T,\rho^{*}}\Big(\|\nabla u\|_{2}^{4}+\|\nabla \omega\|_{2}^{4}\Big)\Big(\|\sqrt{\rho}u_{t}\|_{2}^{2}
+\|\sqrt{\rho}\omega_{t}\|_{2}^{2}\Big),
\end{split}
\end{equation}
and 
\begin{equation}
\begin{split}
I_{5}&\leq\Big|-\int_{\mathbb{T}^{2}}(\sqrt{t}\rho u\cdot\nabla u_{t})\cdot(\sqrt{t}u_{t})dx
-\int_{\mathbb{T}^{2}}(\sqrt{t}\rho u\cdot\nabla \omega_{t})\cdot(\sqrt{t}\omega_{t})dx\Big|\\
&\leq\varepsilon\Big(\|\nabla\sqrt{t} u_{t}\|_{2}^{2}+\|\nabla\sqrt{t} \omega_{t}\|_{2}^{2}\Big)
+C\rho^{*}\|u\|_{\infty}^{2}\Big(\|\sqrt{\rho t}u_{t}\|_{2}^{2}+\|\sqrt{\rho t}\omega_{t}\|_{2}^{2}\Big).
\end{split}
\end{equation}
Therefore,  for some constant $C_{T,\rho^{*}}$ depending only on $\rho^{*}$ and $T$, from \eqref {3.11}, we conclude that
\begin{equation*}
\begin{split}
&\frac{d}{dt}\Big(\|\sqrt{\rho t}u_{t}\|_{2}^{2}+\|\sqrt{\rho t}\omega_{t}\|_{2}^{2}\Big)
+\|\nabla \sqrt{t}u_{t}\|_{2}^{2}+\|\nabla \sqrt{t}\omega_{t}\|_{2}^{2}\\
&\leq C\Big((1+\rho^{*})\|u\|_{\infty}^{2}+\|u\|_{\infty}^{4}\Big)\Big(\|\sqrt{\rho t}u_{t}\|_{2}^{2}+\|\sqrt{\rho t}\omega_{t}\|_{2}^{2}\Big)+C_{T,\rho^{*}}\Big(\|\nabla u\|_{4}^{4}+\|\nabla\omega\|_{4}^{4}+\|\nabla^{2} u\|_{2}^{2}\\
&\quad+\|\nabla^{2}\omega\|_{2}^{2}+\|u\|_{\infty}^{4}(\|\nabla u\|_{2}^{2}+
\|\nabla \omega\|_{2}^{2})+(1+\|\nabla u\|_{2}^{4}+\|\nabla \omega\|_{2}^{4})(\|\sqrt{\rho}u_{t}\|_{2}^{2}
+\|\sqrt{t}\omega_{t}\|_{2}^{2})\Big).
\end{split}
\end{equation*}
Set
\begin{equation*}
\begin{split}
h_1(t)&=C\Big((1+\rho^{*})\|u\|_{\infty}^{2}+\|u\|_{\infty}^{4}\Big)+C_{T,\rho^{*}}\Big(\|\nabla u\|_{4}^{4}+\|\nabla\omega\|_{4}^{4}+\|\nabla^{2} u\|_{2}^{2}+\|\nabla^{2}\omega\|_{2}^{2}\\
&\quad+\|u\|_{\infty}^{4}(\|\nabla u\|_{2}^{2}+
\|\nabla \omega\|_{2}^{2})+(1+\|\nabla u\|_{2}^{4}+\|\nabla \omega\|_{2}^{4})(\|\sqrt{\rho}u_{t}\|_{2}^{2}
+\|\sqrt{t}\omega_{t}\|_{2}^{2})\Big),
\end{split}
\end{equation*}
then, $h_1(t)\in L_{loc}^{1}(\mathbb{R}^{+})$ depending only on $\rho^{*},\|\sqrt{\rho_{0}}u_{0}\|_{2}$, $\|\sqrt{\rho_{0}}\omega_{0}\|_{2}$ and $K_{0}$. Indeed,
from \eqref{3.1}, \eqref{3.2-1} and  the 2-D Gagliardo-Nirenberg interpolation inequality $\|u\|_{\infty}^{4}\leq\|u\|_{2}^{2}\|\nabla^{2}u\|_{2}^{2}$, we get  $(u,\omega)\in L^{4}(\mathbb{R}^{+};L^{\infty})$, and naturally $(u,\omega)\in L^{2}(0,T;L^{\infty}).$ Similarly, we also get $(\nabla u,\nabla\omega)\in L^{4}(0,T;L^{4})$, $(\nabla^{2}u,\nabla^{2}\omega)\in L^{2}(\mathbb{R}^{+};L^{2})$.\\
Obviously, if the solution is smooth with density bounded away from zero, then we have
 \begin{equation*}\label{3.20-A}\lim_{t\rightarrow 0^{+}}\int_{\mathbb{T}^{2}}\rho t(|u_{t}|^{2}+\omega_{t}|^{2})dx=0.\end{equation*}
Thus,  integrating with respect to time from $0$ to $t$ for the following inequality
\begin{equation*}
\begin{split}\label{3.20}
\frac{d}{dt}\Big(\|\sqrt{\rho t}u_{t}\|_{2}^{2}&+\|\sqrt{\rho t}\omega_{t}\|_{2}^{2}
+\int_{0}^{t}\tau\|\nabla u_{t}\|_{2}^{2}d\tau+\int_{0}^{t}\tau\|\nabla\omega_{t}\|_{2}^{2}d\tau\Big)\\
&\leq h_1(t)\Big(1+\|\sqrt{\rho t}u_{t}\|_{2}^{2}+\|\sqrt{\rho t}\omega_{t}\|_{2}^{2}\Big),
\end{split}
\end{equation*}
 we conclude that  \eqref{3.21} holds for $t\geq0$.

As in the process starting from time $t_{0}$, we also have the following lemma.
\begin{Lemma}\label{3.3}
Assume $d=2$ and that the solution is smooth with no vacuum. Then for all $t_{0},T\geq0,$ we have
\begin{align}\label{3.22}
\sup_{t_{0}\leq t\leq t_{0}+T}\int_{\mathbb{T}^{2}}\rho(t-t_{0})\big(|u_{t}|^{2}
+|\omega_{t}|^{2}\big)dx
+\int_{t_{0}}^{t_{0}+T}\int_{\mathbb{T}^{2}}(t-t_{0})\Big(|\nabla u_{t}|^{2}+|\nabla\omega_{t}|^{2}\Big)dxdt
\leq c(T)
\end{align}
with $c(T)$ going to zero as $T\rightarrow 0.$
\end{Lemma}
Furthermore, denoting by $\overline{(u_{t})}$ the average of $u_{t}$,
we have
$$\int_{\mathbb{T}^{2}}\rho u_{t}dx=M\overline{(u_{t})}+\int_{\mathbb{T}^{2}}\rho\big(u_{t}-\overline{(u_{t})}\big)dx.$$
Thus,
$$M|\overline{(u_{t})}|\leq\|\rho\|_{2}\|\nabla u_{t}\|_{2}+M^{\frac{1}{2}}\|\sqrt{\rho} u_{t}\|_{2}.$$
Similarly,
$$M|\overline{(\omega_{t})}|\leq\|\rho\|_{2}\|\nabla \omega_{t}\|_{2}+M^{\frac{1}{2}}\|\sqrt{\rho} \omega_{t}\|_{2}.$$
Adding the above two inequalities yields that
\begin{equation*}M\big(|\overline{(u_{t})}|+|\overline{(\omega_{t})}|\big)\leq\|\rho\|_{2}\Big(\|\nabla u_{t}\|_{2}+\|\nabla \omega_{t}\|_{2}\Big)+M^{\frac{1}{2}}\Big(\|\sqrt{\rho} u_{t}\|_{2}+\|\sqrt{\rho} \omega_{t}\|_{2}\Big).
\end{equation*}
Since $\|\rho\|_{2}$ and $M$ are time independent, we get by Sobolev embedding,   
\begin{equation*}
\begin{split}
\|u_{t}\|_{p}+\|\omega_{t}\|_{p}
&\leq\|u_{t}-\overline{(u_{t})}\|_{p}+|\overline{(u_{t})}|+
\|\omega_{t}-\overline{(\omega_{t})}\|_{p}+|\overline{(\omega_{t})}|\\
&\leq \Big(C_{p}+\frac{\|\rho_{0}\|_{2}}{M}\Big)\Big(\|\nabla u_{t}\|_{2}+\|\nabla \omega_{t}\|_{2}\Big)
+\frac{1}{M^{1/2}}\Big(\|\sqrt{\rho} u_{t}\|_{2}+\|\sqrt{\rho} \omega_{t}\|_{2}\Big),
\end{split}
\end{equation*}
which implies that, for all $p<\infty$,
\begin{equation}
\begin{split}\label{3.23}
&\|\sqrt{t}u_{t}\|_{L^{2}(0,T;L^{p})}+\|\sqrt{t}\omega_{t}\|_{L^{2}(0,T;L^{p})}\\
&\leq\Big(C_{p}+\frac{\|\rho_{0}\|_{2}}{M}\Big)\Big(\|\sqrt{t}\nabla u_{t}\|_{L^{2}(0,T;L^{2})}+\|\sqrt{t}\nabla \omega_{t}\|_{L^{2}(0,T;L^{2})}\Big)\\
&\quad+\frac{1}{M^{1/2}}\Big(\|\sqrt{\rho t} u_{t}\|_{L^{2}(0,T;L^{2})}+\|\sqrt{\rho t}\omega_{t}\|_{L^{2}(0,T;L^{2})}\Big).
\end{split}
\end{equation}
According to  \eqref{3.21}, we deduce that for  $p<\infty$,
\begin{equation}
\label{2.31-1}\|(\sqrt{t}u_{t},\sqrt{t}\omega_{t})\|_{L^{2}(0,T;L^{p})}\leq c(T)\quad \text{with}\quad c(T)\rightarrow 0 \quad \text{for}\quad T\rightarrow 0.\end{equation}

In order to obtain some  strong sense  convergence of the approximate sequences $(u^{\delta}, \omega^{\delta})$, as in \cite{DM},  we also  need the following control on the regularity of $u,\omega$ with respect to the time variable.
\begin{Lemma}\label{4.2}
Let $p\in[1,\infty]$ and $u,\omega$ satisfy $u,\omega \in L^{2}(0,T;L^{p})$ and $\sqrt{t}u_{t}, \sqrt{t}\omega_{t}\in L^{2}(0,T;L^{p})$. Then $u,\omega\in H^{\frac{1}{2}-\alpha}(0,T;L^{p})$ for all $\alpha\in(0,1/2)$ and
\begin{align}\label{4.6}
\|u,\omega\|_{H^{\frac{1}{2}-\alpha}}^{2}\leq\|u,\omega\|_{L^{2}(0,T;L^{p})}^{2}
+C_{\alpha,T}\|\sqrt{t}u_{t}, \sqrt{t}\omega_{t}\|_{L^{2}(0,T;L^{p})}^{2},
\end{align}
with $C_{\alpha,T}$ depending only on $\alpha$ and on $T.$
\end{Lemma}
\par
\section{Existence of solution  and weighted energy method in 3-D }
\subsection{Existence of solution in 3-D}
Similar to 2-D, here,  we only present  the \emph{a priori} estimates for smooth enough solutions $(\rho,u,\omega)$ of  system \eqref{eq:NA} in what follows.
\begin{Proposition}
 \label{Pro:4.1}    Under the assumptions of Theorem \ref{th:main2}.
 Let $(\rho,u,\omega)$ be a smooth enough solution of system \eqref{eq:NA} satisfying \eqref{3.1-1}. There exist a universal positive constant $C$  and  $T>0$ such that
\begin{align}\label{4.3-11111}
T\leq\frac{C}{(\rho^{*})^{3}C_0 K_{0}^{3}}.
\end{align}
Then, for all $t\in [0,T)$, we have
\begin{equation}\label{3.2}
\begin{split}
&\|\nabla u\|_{2}^{2}+\|\nabla \omega\|_{2}^{2}+\int_{0}^{t}\Big(\|\sqrt{\rho}u_{t}\|_{2}^{2}+\|\sqrt{\rho}\omega_{t}\|_{2}^{2}+\|\nabla^{2}u\|_{2}^{2}+\|\nabla^{2}\omega\|_{2}^{2}+\|\nabla P\|^{2}_{2}\Big)d\tau\\
&\leq CK_0.
\end{split}
\end{equation}
Furthermore, if \eqref{2.2} is  satisfied then \eqref{3.2} holds true
 for all $t\in[0,\infty)$.
At last, inequality \eqref{3.3A} holds true for all $p\in[1,6].$
\end{Proposition}
\noindent{\bf Proof.}\ Taking the $L^2$-scalar product of the
first equation of  system \eqref{eq:NA} with $u_{t}$  and
the third equation with $\omega_{t}$ respectively,  we get by  $\rho_{t}=-u\cdot\nabla\rho$ that
\begin{equation*}
\int_{\mathbb{T}^{2}}\rho| u_{t}|^{2}dx+\frac{\nu}{2}\frac{d}{dt}\int_{\mathbb{T}^{2}}|\nabla u|^{2}dx
=2\chi\int_{\mathbb{T}^{2}}\Cl \omega\cdot u_{t}dx-\int_{\mathbb{T}^{2}}(\rho u\cdot\nabla u)\cdot u_{t}dx,
\end{equation*}
and
\begin{equation*}
\begin{split}
\int_{\mathbb{T}^{2}}\rho| \omega_{t}|^{2}dx+&\frac{\gamma}{2}\frac{d}{dt}\int_{\mathbb{T}^{2}}|\nabla \omega|^{2}dx
+\frac{\kappa}{2}\frac{d}{dt}\int_{\mathbb{T}^{2}}|\Dv \omega|^{2}dx
+2\chi\frac{d}{dt}\int_{\mathbb{T}^{2}}|\omega|^{2}dx\\
&=2\chi\int_{\mathbb{T}^{2}}\Cl u\cdot \omega_{t}dx-\int_{\mathbb{T}^{2}}(\rho u\cdot\nabla \omega)\cdot \omega_{t}dx.
\end{split}
\end{equation*}
Adding the two identities above and using H\"{o}lder's and Young's inequalities and the fact $\Cl(\Cl u)=-\Delta u \ (\text{for} \ \Dv u=0)$ yield that 
\begin{equation*}
\begin{split}
&\|\sqrt{\rho}u_{t}\|_{2}^{2}+\|\sqrt{\rho}\omega_{t}\|_{2}^{2}
+\frac{1}{2}\frac{d}{dt}\Big(\mu\|\nabla u\|_{2}^{2}+\gamma\|\nabla \omega\|_{2}^{2}+\kappa\|\Dv \omega\|_{2}^{2}+\chi\|\Cl  u-2\omega\|_{2}^{2}\Big)\\
&=-\int_{\mathbb{T}^{2}}(\rho u\cdot\nabla u)\cdot u_{t}dx-\int_{\mathbb{T}^{2}}(\rho u\cdot\nabla \omega)\cdot \omega_{t}dx\\
&\leq\frac{1}{2}\int_{\mathbb{T}^{2}}\rho|u_{t}|^{2}dx+\frac{1}{2}\int_{\mathbb{T}^{2}}\rho|\omega_{t}|^{2}dx
+\frac{1}{2}\int_{\mathbb{T}^{2}}\rho|u\cdot\nabla u|^{2}dx
+\frac{1}{2}\int_{\mathbb{T}^{2}}\rho|u\cdot\nabla \omega|^{2}dx,
\end{split}
\end{equation*}
which gives that
\begin{equation}\label{3.5}
\begin{split}
&\|\sqrt{\rho}u_{t}\|_{2}^{2}+\|\sqrt{\rho}\omega_{t}\|_{2}^{2}
+\frac{d}{dt}\Big(\mu\|\nabla u\|_{2}^{2}+\gamma\|\nabla \omega\|_{2}^{2}+\kappa\|\Dv \omega\|_{2}^{2}+\chi\|\Cl  u-2\omega\|_{2}^{2}\Big)\\
&\leq\int_{\mathbb{T}^{2}}\rho|u\cdot\nabla u|^{2}dx
+\int_{\mathbb{T}^{2}}\rho|u\cdot\nabla \omega|^{2}dx.
\end{split}
\end{equation}
Furthermore, we also get from integrating with respect to time over $[0,t]$  on the both sides of \eqref{3.5} that  
\begin{equation}\label{3.5-1}
\begin{split}
&\|\nabla u\|_{2}^{2}+\|\nabla \omega\|_{2}^{2}+\int_{0}^{t}\Big(\|\sqrt{\rho}u_{t}\|_{2}^{2}+\|\sqrt{\rho}\omega_{t}\|_{2}^{2}\Big)d\tau\\
&\leq CK_0+C\Big(\int_{0}^{t}\int_{\mathbb{T}^{2}}\rho|u\cdot\nabla u|^{2}dxd\tau
+\int_{0}^{t}\int_{\mathbb{T}^{2}}\rho|u\cdot\nabla \omega|^{2}dxd\tau\Big),
\end{split}
\end{equation}
where $K_0$ is given by \eqref{2.2-2-1}.\\
From \eqref{3.5-1}, \eqref{3.6-1} and \eqref{3.6-2}, we finally conclude that
\begin{equation}\label{3.6-3}
\begin{split}
&\|\nabla u\|_{2}^{2}+\|\nabla \omega\|_{2}^{2}+\int_{0}^{t}\Big(\|\sqrt{\rho}u_{t}\|_{2}^{2}+\|\sqrt{\rho}\omega_{t}\|_{2}^{2}+\|\nabla^{2}u\|_{2}^{2}+\|\nabla^{2}\omega\|_{2}^{2}+\|\nabla P\|^{2}_{2}\Big)d\tau\\
&\leq CK_0+C\Big(\int_{0}^{t}\int_{\mathbb{T}^{2}}\rho|u\cdot\nabla u|^{2}dxd\tau
+\int_{0}^{t}\int_{\mathbb{T}^{2}}\rho|u\cdot\nabla \omega|^{2}dxd\tau\Big).
\end{split}
\end{equation}
In what follows, we will bound the last term on the right-hand side of the inequality above. To this end,  it follows from   H\"{o}lder's and  Young's inequalities and Sobolev embedding $\dot{H}^{1}(\mathbb{T}^{3})\hookrightarrow L^{6}(\mathbb{T}^{3})$  that 
\begin{equation*}\label{4.4-10}
\begin{split}
\int_{\mathbb{T}^{3}}\rho|u\cdot\nabla u|^{2}dx
&\leq (\rho^{*})^{\frac{1}{2}}\|\rho^{\frac{1}{4}} u\|_{4}^{2}\|\nabla u\|_{4}^{2}\\
&\leq (\rho^{*})^{\frac{3}{4}}\|\sqrt{\rho} u\|_{2}^{\frac{1}{2}}\| u\|_{6}^{\frac{3}{2}}\|\nabla u\|_{2}^{\frac{1}{2}}\|\nabla u\|_{6}^{\frac{3}{2}}\\
&\leq (\rho^{*})^{\frac{3}{4}}\|\sqrt{\rho} u\|_{2}^{\frac{1}{2}}\|\nabla u\|_{2}^{\frac{3}{2}}\|\nabla u\|_{2}^{\frac{1}{2}}\|\nabla^{2} u\|_{2}^{\frac{3}{2}}\\
&\leq\varepsilon\|\nabla^{2} u\|_{2}^{2}+C(\rho^{*})^{3}\|\sqrt{\rho} u\|_{2}^{2}\|\nabla u\|_{2}^{8},
\end{split}
\end{equation*}
and
\begin{equation*}\label{4.4-11}
\begin{split}
\int_{\mathbb{T}^{3}}\rho|u\cdot\nabla \omega|^{2}dx
&\leq (\rho^{*})^{\frac{1}{2}}\|\rho^{\frac{1}{4}} u\|_{4}^{2}\|\nabla \omega\|_{4}^{2}\\
&\leq (\rho^{*})^{\frac{3}{4}}\|\sqrt{\rho} u\|_{2}^{\frac{1}{2}}\| u\|_{6}^{\frac{3}{2}}\|\nabla \omega\|_{2}^{\frac{1}{2}}\|\nabla \omega\|_{6}^{\frac{3}{2}}\\
&\leq (\rho^{*})^{\frac{3}{4}}\|\sqrt{\rho} u\|_{2}^{\frac{1}{2}}\|\nabla u\|_{2}^{\frac{3}{2}}\|\nabla \omega\|_{2}^{\frac{1}{2}}\|\nabla^{2} \omega\|_{2}^{\frac{3}{2}}\\
&\leq\varepsilon\|\nabla^{2} \omega\|_{2}^{2}+C(\rho^{*})^{3}\|\sqrt{\rho} u\|_{2}^{2}\|\nabla u\|_{2}^{6}\|\nabla \omega\|_{2}^{2},
\end{split}
\end{equation*}
from which, together with  \eqref{3.1} and \eqref{3.6-3} implies that 
\begin{equation}\label{4.9-1}
\begin{split}
\|\nabla u\|_{2}^{2}+\|\nabla\omega\|_{2}^{2}+
&\int_{0}^{t}\Big(\|\sqrt{\rho}u_{t}\|_{2}^{2}+\|\sqrt{\rho}\omega_{t}\|_{2}^{2}+\|\nabla^{2}u\|_{2}^{2}+\|\nabla^{2}\omega\|_{2}^{2}+\|\nabla P\|^{2}_{2}\Big)d\tau\\
&\leq CK_0+C\int_{0}^{t}(\rho^{*})^{3}\|\sqrt{\rho} u\|_{2}^{2}\|\nabla u\|_{2}^{6}\Big(\|\nabla u\|_{2}^{2}+\|\nabla \omega\|_{2}^{2}\Big)d\tau.
\end{split}
\end{equation}
That is  
\begin{equation}\label{4.9-10}
\begin{split}
X(t)\leq CK_0+C\int_{0}^{t}f_{1}(\tau)X^{3}(\tau)d\tau
\end{split}
\end{equation}
with $f_{1}(t)=(\rho^{*})^{3}\|\sqrt{\rho} u\|_{2}^{2}(\|\nabla u\|_{2}^{2}+\|\nabla \omega\|_{2}^{2})$.\\
Setting $g_1(t)=\int_{0}^{t}f_{1}(\tau)X^{3}(\tau)d\tau$,   from \eqref{4.9-10}, we get
$X(t)\leq CK_0 +Cg_1(t).$  Thus,
\begin{equation*}
\begin{split}\frac{d}{dt}g_1(t)&= f_1(t)X^{3}(t)\\
&\leq f_1(t) \Big(CK_0 +Cg_1(t)\Big)^{3}.
\end{split}
\end{equation*}
Hence, whenever $T$ satisfies  $2CK_0^{2}\int_{0}^{T}f_1(\tau)d\tau\leq\frac{1}{2}$, we obtain
\begin{equation*}
\Big(CK_0 +Cg_1(t)\Big)^{2}\leq \frac{CK_0^{2}}{1-2CK_0^{2}\int_{0}^{t}f_1(\tau)d\tau}\quad \text{for }\  t\in [0,T],
\end{equation*}
which ensures that
\begin{equation}\label{3.10-1-A} X^{2}(t)\leq\frac{CK_0^{2}}{1-2CK_0^{2}\int_{0}^{t}f_1(\tau)d\tau}.\end{equation}
Thus, we conclude that \eqref{3.2} holds for   $t\in [0,T].$
Furthermore, according \eqref{3.1},  we have
$$\int_{0}^{T}f_1(t)dt\leq\big(\rho^{*})^{3}C_0T\sup_{t\in[0,T]}X(t),$$
which implies  that \eqref{4.3-11111} holds.

On the other hand, if $\varepsilon_0$ is small enough in \eqref{2.2}, from  \eqref{3.1} and \eqref{3.10-1-A}, for  $t\in[0,\infty)$,  we also have
$$2CK_0^{2}\int_{0}^{t}f_1(\tau)d\tau\leq\frac{1}{2}, $$
which yields that
\begin{equation*}
\begin{split}
X(t)\leq CK_0\quad \text{for }\quad t\in[0,\infty).
\end{split}
\end{equation*}
The proof of the last part of the theorem is similar to 2-D case. The only difference is that Sobolev embedding $H^{1}(\mathbb{T}^{3})\hookrightarrow L^{p}(\mathbb{T}^{3})$ holds true only for $p\leq6$. Here, we omit it.
\subsection{Weighted energy method in 3-D}
In 3-D case, our aim is also to obtain bounds $(\sqrt{\rho t}u_{t},\sqrt{\rho t}\omega_{t})$ in $L^{\infty}([0,T];L^{2})$  and $(\sqrt{t}\nabla u_{t},\sqrt{t}\nabla \omega_{t})$ in $L^{2}([0,T];L^{2})$ respectively.
\begin{Lemma}\label{3.3-1-1}
Assume $d=3$ and that the solution is smooth enough of  system \eqref{eq:NA} with no vacuum. Then for all $t\geq0,$ we have
\begin{align}\label{3.8-1-1}
\|\sqrt{\rho t}u_{t}\|_{2}^{2}+\|\sqrt{\rho t}\omega_{t}\|_{2}^{2}
+\int_{0}^{t}\tau\|\nabla u_{t}\|_{2}^{2}d\tau+\int_{0}^{t}\tau\|\nabla\omega_{t}\|_{2}^{2}d\tau
\leq \exp\Big(\int_{0}^{t}h_2(\tau)d\tau\Big)-1
\end{align}
with $h_2\in L^{1}_{loc}(\mathbb{R}^{+})$ depending only on $\rho^{*},\|\sqrt{\rho_{0}}u_{0}\|_{2},
\|\sqrt{\rho_{0}}\omega_{0}\|_{2}$ and $K_{0}.$
\end{Lemma}
\noindent{\bf Proof.} Compared with the proof of lemma \ref{3.3}  for the 2-D case, we here only show  some different parts for $I_2$-$I_5$ in what follows.
For $I_{3}$, we also have \eqref{2.27-A}. The method of processing $I_{31}$ and $I_{34}$ in 2-D is not applicable to 3-D. However,  combining H\"{o}lder's inequality and Sobolev embedding $\dot{H}^{1}(\mathbb{T}^{3})\hookrightarrow L^{6}(\mathbb{T}^{3})$  for some constant $C_{T,\rho^{*}}$ depending only on $T$ and $\rho^{*}$, we have
\begin{equation*}
\begin{split}
I_{31}&\leq\sqrt{\rho^{*} T}\|\sqrt{\rho t}u_{t}\|_{4}\|u\|_{6}\|\nabla u\|_{24/7}^{2}\\
&\leq\sqrt{\rho^{*} T}\|\sqrt{\rho t}u_{t}\|_{2}^{1/4}\|\sqrt{\rho t}u_{t}\|_{6}^{3/4}
\|u\|_{6}\|\nabla u\|_{24/7}^{2}\\
&\leq\varepsilon\|\nabla\sqrt{t}u_{t}\|_{2}^{2}+C_{T,\rho^{*}}\|\sqrt{\rho t}u_{t}\|_{2}^{2/5}
\|\nabla u\|_{24/7}^{16/5}\|\nabla u\|_{2}^{8/5}.
\end{split}
\end{equation*}
Due to
$$\|\nabla u\|_{24/7}^{16/5}\leq C\|\nabla u\|_{2}^{6/5}\|\nabla^{2} u\|_{2}^{2},$$
thus,
\begin{equation*}
\begin{split}
I_{31}&\leq\varepsilon\|\nabla\sqrt{t}u_{t}\|_{2}^{2}+C_{T,\rho^{*}}\|\sqrt{\rho t}u_{t}\|_{2}^{2/5}
\|\nabla u\|_{2}^{14/5}\|\nabla^{2} u\|_{2}^{2}\\
&\leq\varepsilon\|\nabla\sqrt{t}u_{t}\|_{2}^{2}+C_{T,\rho^{*}}\Big(\|\sqrt{\rho t}u_{t}\|_{2}^{2}+
\|\nabla u\|_{2}^{7/2}\Big)\|\nabla^{2} u\|_{2}^{2}.
\end{split}
\end{equation*}
Similarly,
\begin{equation*}
\begin{split}
I_{34}&\leq\sqrt{\rho^{*} T}\|\sqrt{\rho t}\omega_{t}\|_{4}\|u\|_{6}\|\nabla u\|_{24/7}\|\nabla \omega\|_{24/7}\\
&\leq\sqrt{\rho^{*} T}\|\sqrt{\rho t}\omega_{t}\|_{2}^{1/4}\|\sqrt{\rho t}\omega_{t}\|_{6}^{3/4}
\|u\|_{6}\|\nabla u\|_{24/7}\|\nabla \omega\|_{24/7}\\
&\leq\varepsilon\|\nabla\sqrt{t}\omega_{t}\|_{2}^{2}+C_{T,\rho^{*}}\|\sqrt{\rho t}\omega_{t}\|_{2}^{2/5}
\|\nabla u\|_{2}^{8/5}\|\nabla u\|_{24/7}^{8/5}\|\nabla \omega\|_{24/7}^{8/5},
\\&\leq\varepsilon\|\nabla\sqrt{t}\omega_{t}\|_{2}^{2}+C_{T,\rho^{*}}\|\sqrt{\rho t}\omega_{t}\|_{2}^{2/5}
\|\nabla u\|_{2}^{11/5}\|\nabla^{2} u\|_{2}\|\nabla \omega\|_{2}^{3/5}\|\nabla^{2} \omega\|_{2}\\
&\leq\varepsilon\|\nabla\sqrt{t}\omega_{t}\|_{2}^{2}+C_{T,\rho^{*}}\Big(\|\sqrt{\rho t}\omega_{t}\|_{2}^{2}+
\|\nabla u\|_{2}^{22/5}+\|\nabla\omega\|_{2}^{2}\Big)\Big(\|\nabla^{2}u\|_{2}^{2}+\|\nabla^{2} \omega\|_{2}^{2}\Big).
\end{split}
\end{equation*}
Other items are treated exactly the same in 3-D as in 2-D.\\
Therefore,  we get  from \eqref{3.11} for some constant $C_{T,\rho^{*}}$ depending only on $\rho^{*}$ and $T$,
\begin{equation*}
\begin{split}
&\frac{d}{dt}\Big(\|\sqrt{\rho t}u_{t}\|_{2}^{2}+\|\sqrt{\rho t}\omega_{t}\|_{2}^{2}\Big)
+\|\nabla \sqrt{t}u_{t}\|_{2}^{2}+\|\nabla \sqrt{t}\omega_{t}\|_{2}^{2}
\\
&\leq C\Big((1+\rho^{*})\|u\|_{\infty}^{2}+\|u\|_{\infty}^{4}\Big)\Big(\|\sqrt{\rho t}u_{t}\|_{2}^{2}+\|\sqrt{\rho t}\omega_{t}\|_{2}^{2}\Big)\\
&\quad+C_{T,\rho^{*}}\Big(\big(\|\sqrt{\rho t}u_{t}\|_{2}^{2}+\|\sqrt{\rho t}\omega_{t}\|_{2}^{2}\big)
\big(\|\nabla^{2} u\|_{2}^{2}+\|\nabla^{2} \omega\|_{2}^{2}\big)+\|u\|_{\infty}^{4}\big(\|\nabla u\|_{2}^{2}+\|\nabla \omega\|_{2}^{2}\big)\\
&\quad+\big(1+\|\nabla u\|_{2}^{7/2}+\|\nabla u\|_{2}^{22/5}+\|\nabla\omega\|_{2}^{2}\big)
\big(\|\nabla^{2} u\|_{2}^{2}+\|\nabla^{2} \omega\|_{2}^{2}\big)
\\
&\quad
+\big(1+\|\nabla u\|_{2}^{4}+\|\nabla\omega\|_{2}^{4}\big)
(\|\sqrt{\rho}u_{t}\|_{2}^{2}+\|\sqrt{\rho}\omega_{t}\|_{2}^{2})\Big).
\end{split}
\end{equation*}
Set
\begin{equation*}
\begin{split}
h_2(t)&=C\big((1+\rho^{*})\|u\|_{\infty}^{2}+\|u\|_{\infty}^{4}\big)+C_{T,\rho^{*}}\Big(\|u\|_{\infty}^{4}\big(\|\nabla u\|_{2}^{2}+\|\nabla \omega\|_{2}^{2}\big)\\
&\quad+\big(1+\|\nabla u\|_{2}^{7/2}+\|\nabla u\|_{2}^{22/5}+\|\nabla\omega\|_{2}^{2}\big)\big(\|\nabla^{2} u\|_{2}^{2}+\|\nabla^{2} \omega\|_{2}^{2}\big)\\
&\quad+\big(1+\|\nabla u\|_{2}^{4}+\|\nabla\omega\|_{2}^{4}\big)
(\|\sqrt{\rho}u_{t}\|_{2}^{2}+\|\sqrt{\rho}\omega_{t}\|_{2}^{2})\Big),
\end{split}
\end{equation*}
thus, $h_2(t)\in L_{loc}^{1}(\mathbb{R}^{+})$ depending only on $\rho^{*},\|\sqrt{\rho_{0}}u_{0}\|_{2}$, $\|\sqrt{\rho_{0}}\omega_{0}\|_{2}$ and $K_{0}$.
In fact, from \eqref{3.1}, \eqref{3.2} and the 3-D Gagliardo-Nirenberg interpolation inequality $\|v\|_{\infty}^{4}\leq\|\nabla v\|_{2}^{2}\|\nabla^{2}v\|_{2}^{2}$, we get  $(u,\omega)\in L^{4}(\mathbb{R}^{+};L^{\infty})$. Similar to the case of 2-D, we can also get that $(u,\omega)\in L^{2}(0,T;L^{\infty})$, $(\nabla u,\nabla\omega)\in L^{4}(0,T;L^{4})$ and $(\nabla^{2}u,\nabla^{2}\omega)\in L^{2}(\mathbb{R}^{+};L^{2})$. \\
Thus, from
\begin{equation*}
\begin{split}\label{3.20-111}
&\frac{d}{dt}\Big(\|\sqrt{\rho t}u_{t}\|_{2}^{2}+\|\sqrt{\rho t}\omega_{t}\|_{2}^{2}
+\int_{0}^{t}\tau\|\nabla u_{t}\|_{2}^{2}d\tau+\int_{0}^{t}\tau\|\nabla\omega_{t}\|_{2}^{2}d\tau\Big)\\
&\leq h_2(t)\Big(1+\|\sqrt{\rho t}u_{t}\|_{2}^{2}+\|\sqrt{\rho t}\omega_{t}\|_{2}^{2}\Big),
\end{split}
\end{equation*}
and  \begin{equation*}\label{3.20-B}\lim_{t\rightarrow 0^{+}}\int_{\mathbb{T}^{23}}\rho t(|u_{t}|^{2}+\omega_{t}|^{2})dx=0,\end{equation*}
we  concludes that  \eqref{3.8-1-1} holds for $t\geq0$.

In a similar way, we know that Lemma \ref{3.3}  still hold  for 3-D, and \eqref{3.23} holds for all $p<6$ in 3-D.
\par
\section{The proof of uniqueness}
The purpose of this section is to present the proof to the uniqueness part of both
Theorems \ref{th:main1} and \ref{th:main2}.
\subsection{More regularity of the solutions}
In order to prove the uniqueness parts of theorems \ref{th:main1} and \ref{th:main2}, we need more information on the regularity of the solution to system \eqref{eq:NA} obtained in previous section.  Our first goal is  to achieve the bound $(\nabla u,\nabla\omega)$ in $L^{1}(0,T;L^{\infty})$ in terms of the data and of $T$ by  performing  the shift of integrability method. This is given by the following two lemmas in 2-D and 3-D respectively.
\begin{Lemma}\label{3.4}
Assume $d=2$, then $\forall T>0,p\in[2,\infty]$ and $\varepsilon$ small enough, we have
\begin{align}\label{3.24}
\|\nabla^{2}\sqrt{t}u,\nabla^{2}\sqrt{t}\omega,\sqrt{t}\omega\|_{L^{p}(0,T;L^{p^{*}-\varepsilon})}
+\|\nabla\sqrt{t}P\|_{L^{p}(0,T;L^{p^{*}-\varepsilon})}\leq C_{0,T},
\end{align}
where $p^{*}\triangleq\frac{2p}{p-2}$ and $C_{0,T}$ depends only on $\rho^{*},\|\sqrt{\rho_{0}} u_{0}\|_{2},\|\sqrt{\rho_{0}} \omega_{0}\|_{2},K_0,p,\varepsilon$.\\
Furthermore,  $\forall s\in [1,2)$, there exists $\theta>0$ such that
\begin{align}\label{3.25}
\int_{0}^{T}\Big(\|\nabla u\|_{\infty}^{s}+\|\nabla \omega\|_{\infty}^{s}\Big)dt\leq C_{0,T}T^{\theta}.
\end{align}
\end{Lemma}
\noindent{\bf Proof.}\ From $\eqref{eq:NA}_{1}$ and $\eqref{eq:NA}_{3}$, we have
\begin{align}\label{3.26}
\left\{
\begin{aligned}
&-\nu\Delta\sqrt{t} u+\nabla\sqrt{t}P=2\chi\Cl(\sqrt{t}\omega)-\rho\sqrt{t}(u_{t}+u\cdot\nabla u)&\quad \hbox{in} \quad (0,T)\times\mathbb{T}^{2},\\
&\Dv\sqrt{t}u =0&\quad \hbox{in} \quad (0,T)\times\mathbb{T}^{2},\\
&-\gamma\Delta\sqrt{t}\omega-\kappa\nabla\Dv(\sqrt{t}\omega)+4\chi\sqrt{t}\omega
=2\chi \Cl(\sqrt{t} u) -\rho\sqrt{t}(\omega_{t}+u\cdot\nabla \omega)&\quad \hbox{in} \quad (0,T)\times\mathbb{T}^{2}.
\end{aligned}
\right.
\end{align}
Using \eqref{3.1-1} and \eqref{3.21} yields that $$(\rho\sqrt{t}u_{t},\rho\sqrt{t}\omega_{t})\in L^{\infty}(0,T;L^{2}).$$
According to \eqref{3.23}, we  have 
$$  (\rho\sqrt{t}u_{t},\rho\sqrt{t}\omega_{t})\in  L^{2}(0,T;L^{q})\quad \text{for}\quad q<\infty.$$ Therefore,  we get by interpolation inequality,
$$\|\rho\sqrt{t}u_{t}\|_{L^{p}(0,T;L^{r})}\leq
\|\rho\sqrt{t}u_{t}\|_{L^{\infty}(0,T;L^{2})}^{1-\frac{2}{p}}
\|\rho\sqrt{t}u_{t}\|_{L^{2}(0,T;L^{q})}^{\frac{2}{p}},$$
$$\|\rho\sqrt{t}\omega_{t}\|_{L^{p}(0,T;L^{r})}\leq
\|\rho\sqrt{t}\omega_{t}\|_{L^{\infty}(0,T;L^{2})}^{1-\frac{2}{p}}
\|\rho\sqrt{t}\omega_{t}\|_{L^{2}(0,T;L^{q})}^{\frac{2}{p}}$$
with $\frac{1}{r}=\frac{p-2}{2p}+\frac{2}{pq},$  $2\leq r< p^{*}$.\\
Thus,
\begin{align}\label{3.27}
\|\rho\sqrt{t}u_{t},\rho\sqrt{t}\omega_{t}\|_{L^{p}(0,T;L^{r})}\leq C_{0,T}\quad \text{for}\quad p\in[2,\infty],\quad r\in [2,p^{*}).
\end{align}
Similarly, it is known from \eqref{3.2-1} that $(\nabla u,\nabla\omega)$ is bounded at $L^{\infty}(0,T;L^{2})\cap L^{2}(0,T;H^{1})$.
By interpolation inequality,   we get  for $\frac{1}{r}=\frac{p-2}{2p}+\frac{2}{pq},$ $2\leq r< p^{*}$
\begin{equation*}
\begin{split}
\|\nabla u\|_{L^{p}(0,T;L^{r})}&\leq
\|\nabla u\|_{L^{\infty}(0,T;L^{2})}^{1-\frac{2}{p}}
\|\nabla u\|_{L^{2}(0,T;L^{q})}^{\frac{2}{p}}\\
&\leq\|\nabla u\|_{L^{\infty}(0,T;L^{2})}^{1-\frac{2}{p}}
\|\nabla u\|_{L^{2}(0,T;H^{1})}^{\frac{2}{p}},
\end{split}
\end{equation*}
\begin{equation*}
\begin{split}
\|\nabla\omega\|_{L^{p}(0,T;L^{r})}&\leq
\|\nabla\omega\|_{L^{\infty}(0,T;L^{2})}^{1-\frac{2}{p}}
\|\nabla\omega\|_{L^{2}(0,T;L^{q})}^{\frac{2}{p}}\\
&\leq\|\nabla\omega\|_{L^{\infty}(0,T;L^{2})}^{1-\frac{2}{p}}
\|\nabla\omega\|_{L^{2}(0,T;H^{1})}^{\frac{2}{p}},
\end{split}
\end{equation*}
which implies that
\begin{align}\label{3.28}
&\|\nabla u,\nabla\omega\|_{L^{p}(0,T;L^{r})}\leq C_{0,T} \quad \text{for }\quad p\geq 2,\quad r<p^{*},
\end{align}
and then
\begin{align}\label{3.28-A}\|\Cl(\sqrt{t}u),\Cl(\sqrt{t}\omega)\|_{L^{p}(0,T;L^{r})}\leq C_{0,T} \quad  \text{for}\quad p\geq 2,\quad r<p^{*}.\end{align}
As obvious, since $u,\omega$ (and thus $\sqrt{t}\rho u,\sqrt{t}\rho\omega$) is bounded in all spaces
$L^{q}(0,T;L^{r})$ (except $q=r=\infty$), we conclude that
\begin{align}\label{3.29}
\|\sqrt{t}\rho u\cdot\nabla u,\sqrt{t}\rho u\cdot\nabla\omega\|_{L^{p}(0,T;L^{r})}\leq C_{0,T}\quad  \text{for }\quad p\in[2,\infty],\quad r\in [2,p^{*}).
\end{align}
Applying the maximal regularity estimate for the Stokes  equations and  the standard  estimate for elliptic equations for \eqref{3.26} yields that
\begin{align}\label{3.30}
\|\nabla^{2}\sqrt{t}u,\nabla^{2}\sqrt{t}\omega,\nabla\sqrt{t}P\|_{L^{p}(0,T;L^{r})}\leq C_{0,T}\quad   \text{for}\quad p\in[2,\infty], \quad r\in [2,p^{*}).
\end{align}
Fix $p\in[2,\infty)$ so that $ps<2(p-s)$ and $1\leq s<2$,  which means that $(\int_{0}^{T}t^{-\frac{ps}{2p-2s}}dt)^{\frac{1}{s}-\frac{1}{p}}\leq C_{0,T}$.  Taking $r\in(2,p^{*})$ such that the embedding $W_{r}^{1}\hookrightarrow L^{\infty}$,  we obtain  from \eqref{3.30}
\begin{equation*}
\begin{split}
&\Big(\int_{0}^{T}\|\nabla u\|_{\infty}^{s}dt\Big)^{\frac{1}{s}}+\Big(\int_{0}^{T}\|\nabla \omega\|_{\infty}^{s})dt\Big)^{\frac{1}{s}}\\
&\leq C\Big(\int_{0}^{T}\big(t^{-1/2}\|\sqrt{t}\nabla u\|_{W_{r}^{1}}\big)^{s}dt\Big)^{\frac{1}{s}}
+C\Big(\int_{0}^{T}\big(t^{-1/2}\|\sqrt{t}\nabla \omega\|_{W_{r}^{1}}\big)^{s}dt\Big)^{\frac{1}{s}}\\
&\leq C\Big(\int_{0}^{T}t^{-\frac{ps}{2p-2s}}dt\Big)^{\frac{1}{s}-\frac{1}{p}}
\Big(\|\nabla\sqrt{t}u\|_{L_{p}(0,T;W_{r}^{1})}+\|\nabla\sqrt{t}\omega\|_{L_{p}(0,T;W_{r}^{1})}\Big)\\
&\leq C_{0,T}T^{\frac{2p-2s-ps}{2ps}},
\end{split}
\end{equation*}
which yields \eqref{3.25}.
\begin{Lemma}\label{4.3} Assume $d=3$, then for all $\ T>0$, $p\in[2,\infty]$, we have
\begin{align}\label{4.7}
\|\nabla^{2}\sqrt{t}u,\nabla^{2}\sqrt{t}\omega,\sqrt{t}\omega\|_{L^{p}(0,T;L^{r})}
+\|\nabla\sqrt{t}P\|_{L^{p}(0,T;L^{r})}\leq C_{0,T},
\quad \text{for} \quad 2\leq r\leq\frac{6p}{3p-4},
\end{align}
where $C_{0,T}$ depends only on $\rho^{*},\|\sqrt{\rho_{0}} u_{0}\|_{2},\|\sqrt{\rho_{0}} \omega_{0}\|_{2},K_{0},p.$\\
Furthermore, for  $s\in [1,\frac{4}{3})$, then for some $\theta>0$,  we have
\begin{align}\label{4.8}
\int_{0}^{T}\Big(\|\nabla u\|_{\infty}^{s}+\|\nabla \omega\|_{\infty}^{s}\Big)dt\leq C_{0,T}T^{\theta}.
\end{align}
\end{Lemma}
\noindent{\bf Proof.}\  From \eqref{3.1-1} and \eqref{3.8-1-1}, we get
$$\rho\sqrt{t}u_{t},\rho\sqrt{t}\omega_{t} \in L^{\infty}(0,T;L^{2})\cap L^{2}(0,T;H^{1}),$$
and from $\dot{H}^{1}(\mathbb{T}^{3}) \hookrightarrow L^{6}(\mathbb{T}^{3})$, we have
$$\rho\sqrt{t}u_{t},\rho\sqrt{t}\omega_{t}\in L^{\infty}(0,T;L^{2})\cap L^{2}(0,T;L^{q})\quad \text{with}\quad q\leq6.$$
It then follows from  interpolation inequality that
$$\|\rho\sqrt{t}u_{t}\|_{L^{p}(0,T;L^{r})}\leq \|\rho\sqrt{t}u_{t}\|_{L^{\infty}(0,T;L^{2})}^{1-\frac{2}{p}}
\|\rho\sqrt{t}u_{t}\|_{L^{2}(0,T;L^{q})}^{\frac{2}{p}},$$ and
$$\|\rho\sqrt{t}\omega_{t}\|_{L^{p}(0,T;L^{r})}\leq \|\rho\sqrt{t}\omega_{t}\|_{L^{\infty}(0,T;L^{2})}^{1-\frac{2}{p}}
\|\rho\sqrt{t}\omega_{t}\|_{L^{2}(0,T;L^{q})}^{\frac{2}{p}},$$
with $\frac{1}{r}=\frac{p-2}{2p}+\frac{2}{pq}$. Here, when $q$ takes 6,  then $r$ may take the maximum value of $\frac{6p}{3p-4}$.\\ Thus, we readily get
\begin{align}\label{4.9}
\|\rho\sqrt{t}u_{t},\rho\sqrt{t}\omega_{t}\|_{L^{p}(0,T;L^{r})}\leq C_{0,T}\quad \text{for}\quad p\in[2,\infty],\quad r\in [2,\frac{6p}{3p-4}].
\end{align}
According to \eqref{3.2}, we have $\nabla u,\nabla\omega\in L^{\infty}(0,T;L^{2})\cap L^{2}(0,T;H^{1}),$ and from $\dot{H}^{1}(\mathbb{T}^{3})\hookrightarrow L^{6}(\mathbb{T}^{3})$, we get $\nabla u,\nabla\omega\in L^{\infty}(0,T;L^{2})\cap L^{2}(0,T;L^{q})$ with $q\leq6.$  By interpolation inequality, we obtain 
$$\|\nabla u\|_{L^{p}(0,T;L^{r})}\leq
\|\nabla u\|_{L^{\infty}(0,T;L^{2})}^{1-\frac{2}{p}}
\|\nabla u\|_{L^{2}(0,T;L^{q})}^{\frac{2}{p}},$$
$$\|\nabla\omega\|_{L^{p}(0,T;L^{r})}\leq
\|\nabla\omega\|_{L^{\infty}(0,T;L^{2})}^{1-\frac{2}{p}}
\|\nabla\omega\|_{L^{2}(0,T;L^{q})}^{\frac{2}{p}}$$
with $\frac{1}{r}=\frac{p-2}{2p}+\frac{2}{pq}.$\\ Then
\begin{align*}
\|\nabla u,\nabla \omega\|_{L^{p}(0,T;L^{r})}\leq C_{0,T}\quad  \text{for }\quad p\in[2,\infty],\quad r\in \big[2,\frac{6p}{3p-4}\big],
\end{align*}
which means that
\begin{equation}\label{3.12-11}
\nabla u,\nabla\omega \in L^{4}(0,T;L^{3}),
\end{equation}
and
\begin{equation*}
\|\Cl(\sqrt{t}u),\Cl(\sqrt{t}\omega)\|_{L^{p}(0,T;L^{r})}\leq C_{0,T}, \quad  \text{for all}\quad p\in[2,\infty],\quad r\in [2,\frac{6p}{3p-4}].
\end{equation*}
On the other hand, using Gagliardo-Nirenberg interpolation inequality $\|v\|_{L^{\infty}}^{4}
\leq C\|\nabla v\|_{L^{2}}^{2}\|\nabla^{2} v\|_{L^{2}}^{2}$ leads to
$$\|u\|_{L^{4}(0,T;L^{\infty})}\leq
\|\nabla u\|_{L^{\infty}(0,T;L^{2})}^{\frac{1}{2}}
\|\nabla^{2} u\|_{L^{2}(0,T;L^{2})}^{\frac{1}{2}},$$
and
$$\|\omega\|_{L^{4}(0,T;L^{\infty})}\leq
\|\nabla \omega\|_{L^{\infty}(0,T;L^{2})}^{\frac{1}{2}}
\|\nabla^{2} \omega\|_{L^{2}(0,T;L^{2})}^{\frac{1}{2}}.$$
Thanks to \eqref{3.1-1} and \eqref{3.2}, we conclude that
\begin{equation}\label{3.12-12}\sqrt{t}\rho u, \sqrt{t}\rho \omega\in L^{4}(0,T;L^{\infty}).\end{equation}
Using  H\"{o}lder's inequality,  and combining with \eqref{3.12-11} and \eqref{3.12-12},  we get
$$\sqrt{t}\rho u\cdot\nabla u,\sqrt{t}\rho u\cdot\nabla\omega\in L^{2}(0,T;L^{3}).$$
Similarly,
$$\sqrt{t}\rho u, \sqrt{t}\rho \omega\in L^{\infty}(0,T;L^{6})\quad \text{and} \quad
\nabla u,\nabla\omega\in L^{\infty}(0,T;L^{2}),$$
which implies  that $$\sqrt{t}\rho u\cdot\nabla u,\sqrt{t}\rho u\cdot\nabla\omega\in L^{\infty}(0,T;L^{3/2}).$$
It follows from  interpolating inequality  and H\"{o}lder's inequality that
\begin{equation*}
\|\sqrt{t}\rho u\cdot\nabla u\|_{L^{p}(0,T;L^{r})}\leq\|\sqrt{t}\rho u\cdot\nabla u\|
_{L^{2}(0,T;L^{3})}^{\frac{2}{p}}
\|\sqrt{t}\rho u\cdot\nabla u\|_{L^{\infty}(0,T;L^{3/2})}
^{1-\frac{2}{p}},
\end{equation*}
\begin{equation*}
\|\sqrt{t}\rho u\cdot\nabla\omega\|_{L^{p}(0,T;L^{r})}\leq\|\sqrt{t}\rho u\cdot\nabla\omega\|_{L^{2}(0,T;L^{3})}^{\frac{2}{p}}
\|\sqrt{t}\rho u\cdot\nabla\omega\|_{L^{\infty}(0,T;L^{3/2})}
^{1-\frac{2}{p}},
\end{equation*}
with $\frac{2}{p}+\frac{3}{r}=2,$  $p\geq2$ .
Using the maximal regularity estimate for the Stokes  equations and  the standard  estimate for elliptic equations for \eqref{3.26}  yields that
\begin{align}\label{4.10}
\begin{aligned}
&\|\nabla^{2}\sqrt{t}u,\nabla^{2}\sqrt{t}\omega,\sqrt{t}\omega\|_{L^{p}(0,T;L^{r})}
+\|\nabla\sqrt{t}P\|_{L^{p}(0,T;L^{r})}\leq C_{0,T}\quad \text{for} \quad p\geq2\quad \text{and}\quad \frac{2}{p}+\frac{3}{r}=2.
\end{aligned}
\end{align}
Furthermore,  using the bound for $(\rho u,\rho\omega)$ in $L^{\infty}(0,T;L^{6})$ and the embedding
$W_{r}^{1}(\mathbb{T}^{3})\hookrightarrow L^{q}(\mathbb{T}^{3})$ with $\frac{3}{q}=\frac{3}{r}-1$ if $1\leq r<3$ (which implies that $(\nabla\sqrt{t}u,\nabla\sqrt{t}\omega)$ is bounded in $L^{p}(0,T;L^{q})$ with $\frac{2}{q}+\frac{3}{r}=2$), we get \eqref{4.10} for the full range of indices.
Fix $p\in(2,4)$ such that $ps<2p-2s$ and take $r=\frac{6p}{3p-4}$. Thanks to $W_{r}^{1}\hookrightarrow L^{\infty}$ (because $r>3$ for $2<p<4$), we have
\begin{equation*}
\begin{split}
&\Big(\int_{0}^{T}\|\nabla u\|_{\infty}^{s}dt\Big)^{\frac{1}{s}}+\Big(\int_{0}^{T}\|\nabla \omega\|_{\infty}^{s})dt\Big)^{\frac{1}{s}}\\
&\leq C\Big(\int_{0}^{T}\|\sqrt{t}\nabla u\|_{W_{r}^{1}}^{s}\frac{dt}{\sqrt{t}}\Big)^{\frac{1}{s}}
+C\Big(\int_{0}^{T}\|\sqrt{t}\nabla \omega\|_{W_{r}^{1}}^{s}\frac{dt}{\sqrt{t}}\Big)^{\frac{1}{s}}\\
&\leq C\Big(\int_{0}^{T}t^{-\frac{ps}{2p-2s}}dt\Big)^{\frac{1}{s}-\frac{1}{p}}
\Big(\|\nabla\sqrt{t}u\|_{L_{p}(0,T;W_{r}^{1})}+\|\nabla\sqrt{t}\omega\|_{L_{p}(0,T;W_{r}^{1})}\Big)\\
&\leq C_{0,T}T^{\frac{2p-2s-ps}{2ps}},
\end{split}
\end{equation*}
which concludes that  \eqref{4.8} holds.
\subsection{Lagrangian formulation}
As in \cite{Dan3,BD,PZZ}, we shall prove the uniqueness part of  both
Theorems \ref{th:main1} and \ref{th:main2} using the Lagrangian formulation of system \eqref{eq:NA}. First, we introduce the flow $X:\mathbb{R}^{+}\times\mathbb{T}^{d}\rightarrow\mathbb{T}^{d}$
 of $u$ by $$\partial_{t}X(t,y)=u(t,X(t,y)),\quad X(0,y)=y. $$
 Note that
$$X(t,y)=y+\int_{0}^{t}u(\tau,X(\tau,y))d\tau,$$
and 
$$\nabla_{y}X(t,y)=Id+\int_{0}^{t}\nabla_{y}u(\tau,X(\tau,y))d\tau.$$
In Lagrangian coordinates $(t, y)$,  a solution $(\rho,u,\omega,P)$ to system $\eqref{eq:NA}$ recasts in $(\bar{\rho},\bar{u},\bar{\omega},\bar{P})$ with
\begin{align}\label{5.1}
\begin{aligned}
&\bar{\rho}(t,y)=\rho(t,X(t,y)),\quad \bar{u}(t,y)=u(t,X(t,y)),\\
&\bar{\omega}(t,y)=\omega(t,X(t,y)), \quad \bar{P}(t,y)=P(t,X(t,y)),
\end{aligned}
\end{align}
and the triplet $(\bar{\rho},\bar{u},\bar{\omega},\bar{P})$ thus satisfies
\begin{align}\label{5.3}
\left\{
\begin{aligned}
&\bar{\rho}\bar{u}_{t}-\nu\Delta_{u}\bar{u}+\nabla_{u}\bar{P}=2\chi\Cl_{u}\bar{\omega}
&\quad \text{in} \quad (0,T)\times\mathbb{T}^{d},\\
&\Dv_{u}\bar{u} =0 &\quad \text{in} \quad (0,T)\times\mathbb{T}^{d}, \\
&\bar{\rho}\bar{\omega}_{t}-\gamma\Delta_{u}\bar{\omega}-\kappa\nabla_{u}\Dv_{u}\bar{\omega}
+4\chi\bar{\omega}=2\chi \Cl_{u}\bar{u} &\quad \text{in} \quad (0,T)\times\mathbb{T}^{d},\\
&\bar{\rho}_{t}=0 &\quad \text{in} \quad (0,T)\times\mathbb{T}^{d},\\
&\bar{\rho}(y,0)=\rho_{0}(y), \bar{u}(y,0)=u_{0}(y), \bar{\omega}(y,0)=\omega_{0}(y) &\quad y\in\mathbb{T}^{d},
\end{aligned}
\right.
\end{align}
where operators $\nabla_{u}, \Delta_{u}, \nabla_{u}\Dv_{u}, \Cl_{u} $ and $\Dv_{u}$ correspond to the original operators $\nabla, \Delta, \nabla\Dv, \Cl $ and $\Dv$, respectively, after performing the change to the Lagrangian coordinates.

As pointed out in \cite{Dan3,BD,PZZ}, in our regularity framework, that latter system \eqref{5.3} is equivalent to system
$\eqref{eq:NA}$. Thanks to \eqref{3.25}  and \eqref{4.8}, we can take the time
$T$ to be small enough so that
\begin{align}\label{5.4}
\int_{0}^{t}\|\nabla u\|_{\infty}d\tau\leq\frac{1}{2}.
\end{align}
Set
\begin{align*}
\begin{aligned}
&A=(\nabla X)^{-1}(\text{inverse of deformation tensor}),\\
&J=\text{det}\nabla X(\text{Jacobian determinant}),\\
&a=JA(\text{tranpose of cofactor matrix}).\\
\end{aligned}
\end{align*}
Thus,  in the $(t, y)$-coordinates, operators $\nabla,\Dv,\Cl$ and $\Delta$ translate into
\begin{align}\label{5.2}
\nabla_{u}:={}^{T}A\nabla_{y},\quad \Dv_{u}:=\Dv_{y}(A\cdot),\quad \Cl_{u}:=\nabla_{u}\wedge\cdot,\quad \text{and} \quad\Delta_{u}:=\Dv_{u}
\nabla_{u}.
\end{align}
Moreover, given some matrix $N,$ we  define the divergence operator (acting on vector fields $v$) by the formulation
\begin{equation}\label{5.2-1}
\div_{u}N v=\Dv_{y}(N\cdot v)\eqdefa{}^T\!N:\nabla v,
\end{equation}
where $N:B=\sum_{i,j}N_{ij}B_{ji}$ for $N=(N_{ij})_{1\leq i,j\leq d}$ and $B=(B_{ij})_{1\leq i,j\leq d}$ two
$d\times d$ matrices.\\
Of course, if the condition \eqref{5.4} is fulfilled then  we have
\begin{align}\label{5.5}
A=\Big(Id+(\nabla_{y}X-Id)\Big)^{-1}=\sum_{k=0}^{+\infty}(-1)^{k}\Big(\int_{0}^{t}\nabla_{y}\bar{u}(\tau,\cdot)d\tau\Big)^{k},
\end{align}
which yields that
\begin{align}\label{5.14}
\delta A=\Big(\int_{0}^{t}\nabla\delta ud\tau\Big)\cdot\Big(\sum_{k\geq1}\sum_{0\leq j<k}C_{1}^{j}C_{2}^{k-1-j}\Big)\quad \text{with} \quad C_{i}(t)=\int_{0}^{t}\nabla\bar{u}^{i}d\tau,
\end{align}
where $\delta A\eqdefa A_{2}-A_{1}$ and $\delta u \eqdefa\bar{u}^{2}-\bar{u}^{1}.$

We also make use of the following permutation symbol
\begin{equation*}
\varepsilon_{ijk}=\begin{cases}
1,& \text{even permutation of {1,2,3}},\\
-1,& \text{odd permutation of {1,2,3}},\\
0,& \text{otherwise},
\end{cases}
\end{equation*}
and the basic identity regarding the $i^{th}$ component of the $\Cl$ of a vector field $u$
$$\big(\Cl u\big)_{i}=\varepsilon_{ijk}u^{k}_{,j}.$$
The chain rule shows that
\begin{equation}\label{curl}\big(\Cl u(X)\big)_{i}=\big(\Cl_{u}\bar{u}\big)_{i}:=\varepsilon_{ijk}A_{j}^{s}\bar{u}^{k}_{,s}.
\end{equation}
Here, we also  present  the following   Piola identity,  that is,  the columns of every cofactor matrix are divergence-free and satisfy
\begin{equation}\label{Piola} a_{i ,k}^{k}=0.\end{equation}
Here, it is pointed out that we use the notation $F_{,k}$ to denote $\frac{\partial F}{\partial x_{k}},$ the $k^{th}$-partial derivative of $F$ for $k=1,2,3$, and omit Einstein's summation convention in \eqref{curl} and \eqref{Piola}.
\subsection{The proof of the uniqueness}
Let $(\rho^{1},u^{1},\omega^{1},P^{1})$ and $(\rho^{2},u^{2},\omega^{2},P^{2})$ be two solutions of   system  $\eqref{eq:NA}$ fulfilling the properties of Theorems \ref{th:main1} and \ref{th:main2},  with the same initial data, and denote by $(\bar{\rho}^{1},\bar{u}^{1},\bar{\omega}^{1},\bar{P}^{1})$ and $(\bar{\rho}^{2},\bar{u}^{2},\bar{\omega}^{2},\bar{P}^{2})$  in Lagrangian coordinates. Of course, we have $\bar{\rho}^{1}=\bar{\rho}^{2}=\rho_{0}$, which explains the choice of our approach here. In what follows, we shall use repeatedly the fact that for $i=1,2$, we have
\begin{align}\label{5.6}
\begin{aligned}
&t^{\frac{1}{2}}\nabla \bar{u}^{i}\in L^{2}(0,T;L^{\infty}),\quad t^{\frac{1}{2}}\nabla \bar{\omega}^{i}\in L^{2}(0,T;L^{\infty}),\quad t^{\frac{1}{2}}\nabla \bar{P}^{i}\in L^{2}(0,T;L^{3}),
\quad t^{\frac{1}{2}}\bar{u}_{t}^{i}\in L^{4/3}(0,T;L^{6}),\\
&
\quad\nabla \bar{u}^{i}\in L^{1}(0,T;L^{\infty})\cap L^{2}(0,T;L^{6})\cap L^{4}(0,T;L^{3}),
\quad \bar{u}^{i}\in L^{4}(0,T;L^{\infty}).
\end{aligned}
\end{align}
It should be noted that the former four items in \eqref{5.6} is less than or equal to $c(T)$, where $c(T)$ designates a nonnegative continuous increasing function of $T$, with
$c(0)=0$ and $c(T)\rightarrow0$ when $T\rightarrow0$. For example, in 3-D,  using   Gagliardo-Nirenberg interpolation inequality, \eqref{3.2-1} and Lemma \ref{4.3}, we have
\begin{equation*}
\begin{aligned}
\|t^{\frac{1}{2}}\nabla \bar{u}^{i}\|^{2}_{L^{2}(0,T;L^{\infty})}
&=\int_{0}^{T}t\|{}^{T}A\nabla u^{i}\|_{\infty}^{2}dt\\
&\leq C\int_{0}^{T}t\|\nabla u^{i}\|_{L^{2}}^{\frac{1}{2}}
\|\nabla^{2} u^{i}\|_{6}^{\frac{3}{2}}dt\\
&\leq C\sup_{t\in[0,T]}\|\nabla u^{i}\|_{L^{2}}^{\frac{1}{2}}
\int_{0}^{T}t\|\nabla^{2} u^{i}\|_{6}^{\frac{3}{2}}dt\\
&\leq CT^{\frac{1}{2}}\|\sqrt{t}\nabla^{2} u^{i}\|_{L^{2}(0,T;L^{6})}^{\frac{3}{2}}\\&\leq c(T).
\end{aligned}
\end{equation*}
Denoting $\delta \omega\eqdefa\bar{\omega}^{2}-\bar{\omega}^{1}$, and $\delta P\eqdefa\bar{P}^{2}-\bar{P}^{1}$, we get
\begin{align}\label{5.7}
\left\{
\begin{aligned}
&\rho_{0}\partial_{t}\delta u-\nu\Delta_{u^{1}}\delta u+\nabla_{u^{1}}\delta P
-2\chi\Cl_{u^{1}}\delta\omega=\delta f_{1},\\
&\Dv_{u^{1}}\delta u =(\Dv_{u^{1}}-\Dv_{u^{2}})\bar{u}^{2} , \\
&\rho_{0}\partial_{t}\delta\omega-\gamma\Delta_{u^{1}}\delta\omega
-\kappa\nabla_{u^{1}}\Dv_{u^{1}}\delta\omega+4\chi\delta\omega
-2\chi\Cl_{u^{1}}\delta u=\delta f_{2},\\
&(\delta u,\delta\omega)|_{t=0}=(0,0),
\end{aligned}
\right.
\end{align}
with $\delta f_{1}\eqdefa\nu(\Delta_{u^{2}}-\Delta_{u^{1}})\bar{u}^{2}-(\nabla_{u^{2}}
-\nabla_{u^{1}})\bar{P}^{2}+2\chi(\Cl_{u^{2}}-\Cl_{u^{1}})\bar{\omega}^{2},$\\
\qquad $\delta f_{2}\eqdefa\gamma(\Delta_{u^{2}}-\Delta_{u^{1}})\bar{\omega}^{2}
+\kappa(\nabla_{u^{2}}\Dv_{u^{2}}-\nabla_{u^{1}}\Dv_{u^{1}})\bar{\omega}^{2}
+2\chi(\Cl_{u^{2}}-\Cl_{u^{1}})\bar{u}^{2}.$\\
We claim for sufficiently small $T>0$,
$$\int_{0}^{T}\int_{\mathbb{T}^{d}}\Big(|\delta u(t,y)|^{2}+|\delta \omega(t,y)|^{2}+|\nabla\delta u(t,y)|^{2}+|\nabla\delta \omega(t,y)|^{2}\Big)dydt=0.$$
To prove our claim,  we first decompose $\delta u$ into
\begin{align}\label{5.8}
\begin{aligned}
\delta u=\varphi+\phi,
\end{aligned}
\end{align}
with $\varphi$ is the solution given by Lemma \ref{l:div} to the following problem:
\begin{align}\label{5.9}
\Dv_{u^{1}}\varphi=(\Dv_{u^{1}}-\Dv_{u^{2}})\bar{u}^{2}=\Dv(\delta A\bar{u}^{2}).
\end{align}
Then, \eqref{eq:diveq} and  \eqref{5.5} ensure that there exist two universal positive constants $c$ and $C$ such that if
\begin{align}\label{5.11}
\|\nabla\bar{u}^{1}\|_{L^{1}(0,T;L^{\infty})}+\|\nabla\bar{u}^{1}\|_{L^{2}(0,T;L^{6})}\leq c,
\end{align}
and then the following inequalities hold true:
\begin{equation}\label{5.12}
\begin{aligned}
&\|\varphi\|_{L^{4}(0,T;L^{2})}\leq C\|\delta A\bar{u}^{2}\|_{L^{4}(0,T;L^{2})},\quad
\|\nabla \varphi\|_{L^{2}(0,T;L^{2})}\leq C\|{}^{T}\delta A:\nabla\bar{u}^{2}\|_{L^{2}(0,T;L^{2})}\\
&\text{and}\quad \|\varphi_{t}\|_{L^{4/3}(0,T;L^{3/2})}\leq C\|\delta A\bar{u}^{2}\|_{L^{4}(0,T;L^{2})}+
C\|(\delta A\bar{u}^{2})_{t}\|_{L^{4/3}(0,T;L^{3/2})}.
\end{aligned}
\end{equation}
 Now, let us bound the r.h.s. of \eqref{5.12}. Regarding ${}^{T}\delta A:\nabla\bar{u}^{2}$, it follows from H\"{o}lder's inequality, \eqref{5.11} and \eqref{5.14} that
\begin{align}\label{5.13}
\sup _{t\in[0,T]}\|t^{-1/2}\delta A\|_{2}\leq C\sup _{t\in[0,T]}\|t^{-1/2}\int_{0}^{t}\nabla\delta ud\tau\|_{2}\leq C\|\nabla \delta u\|_{L_{2}(0,T;L_{2})}.
\end{align}
According  to \eqref{5.6} and \eqref{5.13}, we obtain
\begin{equation*}
\begin{split}
\|{}^{T}\delta A:\nabla\bar{u}^{2}\|_{L^{2}(0,T\times\mathbb{T}^{d})}
&\leq\sup _{t\in[0,T]}\|t^{-1/2}\delta A\|_{2}\|t^{1/2}\nabla\bar{u}^{2}\|_{L^{2}(0,T;L^{\infty})}\\
&\leq c(T)\|\nabla \delta u\|_{L^{2}(0,T;L^{2})}.
\end{split}
\end{equation*}
Similarly,
$$\|\delta A\bar{u}^{2}\|_{L^{4}(0,T;L^{2})}\leq\|t^{-1/2}\delta A\|_{L^{\infty}(0,T;L^{2})}\|t^{1/2}\bar{u}^{2}\|_{L^{4}(0,T;L^{\infty})}.$$
Using \eqref{5.6}, \eqref{5.12} and \eqref{5.13} yields that
\begin{align}\label{5.15-0}\|\nabla \varphi\|_{L^{2}(0,T;L^{2})}\leq c(T)\|\nabla \delta u\|_{L^{2}(0,T;L^{2})},\end{align} and
\begin{align}\label{5.15}
\|\varphi\|_{L^{4}(0,T;L^{2})}\leq c(T)\|\nabla \delta u\|_{L^{2}(0,T;L^{2})}.
\end{align}
In order to bound $\varphi_{t}$, it suffices to derive an appropriate estimate in $L^{4/3}(0, T;L^{3/2})$ for
$$(\delta A\bar{u}^{2})_{t}=\delta A\bar{u}_{t}^{2}+(\delta A)_{t}\bar{u}^{2}.$$
Thanks to \eqref{5.6} and \eqref{5.13}, we have
\begin{equation*}
\begin{split}
\|\delta A\bar{u}_{t}^{2}\|_{L^{4/3}(0,T;L^{3/2})}
&\leq\|t^{-1/2}\delta A\|_{L^{\infty}(0,T;L^{2})}
\|t^{1/2}\nabla\bar{u}_{t}^{2}\|_{L^{4/3}(0,T;L^{6})}\\
&\leq c(T)\|\nabla \delta u\|_{L^{2}(0,T;L^{2})}.
\end{split}
\end{equation*}
For the term $(\delta A)_{t}\bar{u}^{2}$, it follows from H\"{o}lder's inequality that
$$\|(\delta A)_{t}\bar{u}^{2}\|_{L^{4/3}(0,T;L^{3/2})}
\leq\|(\delta A)_{t}\|_{L^{2}(0,T\times\mathbb{T}^{d})}
\|\bar{u}^{2}\|_{L^{4}(0,T;L^{6})}.$$
Furthermore, differentiating \eqref{5.14} with respect to $t$ and using \eqref{5.11} for $\bar{u}^{1}$ and $\bar{u}^{2}$ yield that
$$\|(\delta A)_{t}\|_{2}\leq C\Big(\|\nabla\delta u\|_{2}
+\|t^{-1/2}\int_{0}^{t}\nabla\delta ud\tau\|_{2}\big(\|t^{1/2}\nabla \bar{u}^{1}\|_{\infty}+\|t^{1/2}\nabla \bar{u}^{2}\|_{\infty}\big)\Big).$$
This gives 
$$\|(\delta A)_{t}\|_{L_{2}(0,T\times\mathbb{T}^{d})}\leq C\|\nabla \delta u\|_{L^{2}(0,T\times\mathbb{T}^{d})},$$
which together with \eqref{5.6} implies that
$$\|(\delta A)_{t}\bar{u}^{2}\|_{L^{4/3}(0,T;L^{3/2})}\leq c(T)\|\nabla \delta u\|_{L^{2}(0,T\times\mathbb{T}^{d})}.$$
Thus,
\begin{align}\label{5.15-1}\|\varphi_{t}\|_{L^{4/3}(0,T;L^{3/2})}\leq c(T)\|\nabla \delta u\|_{L^{2}(0,T;L^{2})}.\end{align}
Combining with \eqref{5.12}, \eqref{5.15-0}, \eqref{5.15} and \eqref{5.15-1}, we have
\begin{align}\label{5.10}
\|\varphi\|_{L^{4}(0,T;L^{2})}+\|\nabla \varphi\|_{L^{2}(0,T\times\mathbb{T}^{d})}
+\|\varphi_{t}\|_{L^{4/3}(0,T;L^{3/2})}\leq c(T)\|\nabla \delta u\|_{L^{2}(0,T\times\mathbb{T}^{d})}.
\end{align}
Next, let us restate the equations for $(\delta u,\delta\omega,\delta P)$ as the following system for $(\phi,\delta\omega,\delta P)$:
\begin{align}\label{5.16}
\left\{
\begin{aligned}
&\rho_{0}\partial_{t}\phi-\nu\Delta_{u^{1}}\phi+\nabla_{u^{1}}\delta P
=\delta f_{1}-\rho_{0}\partial_{t}\varphi+\nu\Delta_{u^{1}}\varphi+2\chi\Cl_{u^{1}}\delta\omega,\\
&\Dv_{u^{1}}\phi=0, \\
&\rho_{0}\partial_{t}\delta\omega-\gamma\Delta_{u^{1}}\delta\omega
-\kappa\nabla_{u^{1}}\Dv_{u^{1}}\delta\omega+4\chi\delta\omega
=\delta f_{2}+2\chi\Cl_{u^{1}}\varphi+2\chi\Cl_{u^{1}}\phi.
\end{aligned}
\right.
\end{align}
 Due to  $\Dv_{u^{1}}\phi=0$, we have
\begin{equation*}\label{5.17}
\int_{\mathbb{T}^{d}}(\nabla_{u^{1}} \delta P)\cdot \phi dx=-\int_{\mathbb{T}^{d}}\Dv_{u^{1}}\phi \cdot\delta P dx=0.
\end{equation*}
Note that
\begin{equation*}
\begin{aligned}
&2\chi\int_{\mathbb{T}^{d}}\Cl_{u^{1}}\delta\omega\cdot \phi dx
+2\chi\int_{\mathbb{T}^{d}}\Cl_{u^{1}}\phi\cdot\delta\omega dx\\
&=4\chi\int_{\mathbb{T}^{d}}\Cl_{u^{1}}\phi\cdot\delta\omega dx\\
&\leq4\chi\|\nabla_{u^{1}}\phi\|_2\|\delta\omega\|_2\\
&\leq\chi\|\nabla_{u^{1}}\phi\|_2^2+4\chi\|\delta\omega\|_2^2.
\end{aligned}
\end{equation*}
Therefore, taking the $L^2$-scalar product of the first equation to system \eqref{5.16} with $\phi$  and
the third equation with $\delta\omega$ respectively yields that
\begin{equation}
\begin{split}\label{5.20}
\frac{1}{2}\frac{d}{dt}&\int_{\mathbb{T}^{d}}\rho_{0}\big(|\phi|^{2}+|\delta\omega|^{2}\big) dx
+\int_{\mathbb{T}^{d}}\big(\mu|\nabla_{u^{1}}\phi|^{2}+\gamma|\nabla_{u^{1}}\delta\omega|^{2}
+\kappa|\Dv_{u^{1}}\delta\omega|^{2} \big)dx\\
&=-\int_{\mathbb{T}^{d}}\rho_{0}\partial_{t}\varphi\cdot \phi dx+\nu\int_{\mathbb{T}^{d}}\Delta_{u^{1}}\varphi\cdot \phi dx+\nu\int_{\mathbb{T}^{d}}(\Delta_{u^{2}}-\Delta_{u^{1}})\bar{u}^{2}\cdot \phi dx\\
&\quad-\int_{\mathbb{T}^{d}}(\nabla_{u^{2}}-\nabla _{u^{1}})\bar{P}^{2}\cdot \phi dx+2\chi\int_{\mathbb{T}^{d}}(\Cl_{u^{2}}-\Cl _{u^{1}})\bar{\omega}^{2}\cdot \phi dx\\
&\quad+2\chi\int_{\mathbb{T}^{d}}\Cl_{u^{1}}\varphi\cdot\delta\omega dx+\gamma\int_{\mathbb{T}^{d}}\big(\Delta_{u^{2}}-\Delta_{u^{1}}\big)\bar{\omega}^{2}\cdot\delta\omega dx\\
&\quad+\kappa\int_{\mathbb{T}^{d}}\big(\nabla_{u^{2}}\Dv_{u^{2}}-\nabla_{u^{1}}\Dv_{u^{1}}\big)\bar{\omega}^{2}\cdot\delta\omega dx+2\chi\int_{\mathbb{T}^{d}}\big(\Cl_{u^{2}}-\Cl_{u^{1}}\big)\bar{u}^{2}\cdot\delta\omega dx
\\&\triangleq \sum_{k=1}^{9}II_{k}.
\end{split}
\end{equation}
Here and in what follows,  we estimate term by term above. For $II_{1}$,  it follows  from H\"{o}lder's inequality that
$$\int_{0}^{T}II_{1}(t)dt\leq \|\rho_{0}\|_{\infty}^{3/4}\|\varphi_{t}\|_{L^{4/3}(0,T;L^{3/2})}\|\rho_{0}^{1/4}\phi\|_{L^{4}(0,T;L^{3})}.$$
Using H\"{o}lder's inequality and the Sobolev embedding $H^{1}(\mathbb{T}^{d})\hookrightarrow L^{6}(\mathbb{T}^{d})$ yields that
$$\|\rho_{0}^{1/4}\phi\|_{L^{4}(0,T;L^{3})}\leq\|\sqrt{\rho_{0}}\phi\|_{L^{\infty}(0,T;L^{2})}^{1/2}
\|\phi\|_{L^{2}(0,T;L^{6})}^{1/2}\leq C\|\sqrt{\rho_{0}}\phi\|_{L^{\infty}(0,T;L^{2})}^{1/2}
\|\phi\|_{L^{2}(0,T;H^{1})}^{1/2}.$$
Employing Poincar\'e's inequality in  the unit torus $\mathbb{T}^d$ in \cite{BD}:
$\|\phi\|_{H^{1}}\leq C(\|\sqrt{\rho_{0}}\phi\|_{2}+\|\nabla \phi\|_{2})
$ with constant $C$ depending only on $\rho_{0}$,
 and  taking advantage of \eqref{5.14} and \eqref{5.10}, we conclude that
$$\int_{0}^{T}II_{1}(t)dt\leq c(T)\Big(\|\sqrt{\rho_{0}}\phi\|_{L^{\infty}(0,T;L^{2})}+\|\nabla \phi\|_{L^{2}(0,T\times\mathbb{T}^{d})}\Big)^{1/2}\|\sqrt{\rho_{0}}\phi\|_{L^{\infty}(0,T;L^{2})}^{1/2}
\|\nabla \delta u\|_{L^{2}(0,T\times\mathbb{T}^{d})}.$$
For $II_{2}$, it follows  from integrating by parts and using \eqref{5.10} that
\begin{equation*}
\begin{split}
\int_{0}^{T}II_{2}(t)dt
&\leq\nu\int_{0}^{T}\Big|\int_{\mathbb{T}^{d}}\nabla_{u^{1}}\varphi\nabla_{u^{1}} \phi dx\Big|dt\\
&\leq\nu\int_{0}^{T}\int_{\mathbb{T}^{d}}|\nabla_{u^{1}}\varphi||\nabla_{u^{1}}\phi |dxdt\\
&\leq\frac{\nu}{2}\int_{0}^{T}\|\nabla_{u^{1}}\phi \|_{2}^{2}dt
+\frac{\nu}{2}\int_{0}^{T}\|\nabla_{u^{1}}\varphi\|_{2}^{2}dt\\
&\leq \frac{\nu}{2}\int_{0}^{T}\|\nabla_{u^{1}}\phi \|_{2}^{2}dt+
c(T)\int_{0}^{T}\|\nabla \delta u\|_{2}^{2}dt.
\end{split}
\end{equation*}
For $II_{3}$, using \eqref{5.2} and \eqref{5.14}, we obtain
\begin{equation*}
\begin{split}
II_{3}&\leq C\Big|\int_{\mathbb{T}^{d}}\Dv((\delta A{}^{T}A_{2}+A_{1}{}^{T}\delta A)
\nabla\bar{u}^{2})\cdot \phi dx\Big|\\
&\leq C\int_{\mathbb{T}^{d}}|\delta A{}^{T}A_{2}+A_{1}{}^{T}\delta A||\nabla\bar{u}^{2}||\nabla \phi|dx\\
&\leq C\|t^{-1/2}\delta A\|_{2}\|t^{1/2}\nabla\bar{u}^{2}\|_{\infty}
\|\nabla \phi\|_2,
\end{split}
\end{equation*}
which together with \eqref{5.6} and \eqref{5.13} implies that
\begin{equation*}
\begin{split}\label{5.23}
\int_{0}^{T}II_{3}(t)dt
&\leq\|t^{-1/2}\delta A\|_{L^{\infty}(0,T;L^{2})}
\|t^{1/2}\nabla\bar{u}^{2}\|_{L^{2}(0,T;L^{\infty})}
\|\nabla \phi\|_{L^{2}(0,T\times\mathbb{T}^{d})}\\
&\leq c(T)\|\nabla \delta u\|_{L^{2}(0,T;L^{2})}\|\nabla \phi\|_{L^{2}(0,T\times\mathbb{T}^{d})}.
\end{split}
\end{equation*}
For $II_{4}$, using H\"{o}lder's inequality, we obtain
\begin{align*}\label{5.24}
II_{4}(t)\leq\Big|\int_{\mathbb{T}^{d}}\delta A\nabla\bar{P}^{2}\cdot \phi dx\Big|
\leq C\|t^{-1/2}\delta A\|_{2}\|t^{1/2}\nabla\bar{P}^{2}\|_{3}\|\phi\|_{6}.
\end{align*}
It then follows  from \eqref{5.6}, \eqref{5.13} and Sobolev embedding that
\begin{equation*}
\begin{split}
\int_{0}^{T}II_{4}(t)dt
&\leq\|t^{-1/2}\delta A\|_{L^{\infty}(0,T;L^{2})}
\|t^{1/2}\nabla\bar{P}^{2}\|_{L^{2}(0,T;L^{3})}
\|\phi\|_{L^{2}(0,T;H^{1})}\\
&\leq c(T)\|\nabla \delta u\|_{L^{2}(0,T;L^{2})}
\Big(\|\sqrt{\rho_0}\phi\|_{L_{\infty}(0,T;L_{2})}+\|\nabla \phi\|_{L^{2}(0,T\times\mathbb{T}^{d})}\Big).
\end{split}
\end{equation*}
For $II_{5}$, thanks to  \eqref{5.2}, \eqref{5.6} and H\"{o}lder's inequality, we get
\begin{equation*}
\begin{split}
\int_{0}^{T}II_{5}(t)dt
&\leq C\int_{0}^{T}\Big|\int_{\mathbb{T}^{d}}{}^{T}\delta A\nabla\wedge\bar{\omega}^{2}\cdot \phi dx\Big|dt\\
&\leq C\|t^{-\frac{1}{2}}{}^{T}\delta A\|_{L^{\infty}(0,T;L^{2})}
\|t^{\frac{1}{2}}\nabla\bar{\omega}^{2}\|_{L^{2}(0,T;L^{4})}
\|\phi\|_{L^{2}(0,T;L^{4})}\\
&\leq c(T)\|\nabla \delta u\|_{L^{2}(0,T;L^{2})}\|\phi\|_{L^{2}(0,T;H^{1})}\\
&\leq c(T)\|\nabla \delta u\|_{L^{2}(0,T;L^{2})}
\Big(\|\sqrt{\rho_0}\phi\|_{L_{\infty}(0,T;L_{2})}+\|\nabla \phi\|_{L^{2}(0,T\times\mathbb{T}^{d})}\Big).
\end{split}
\end{equation*}
For $II_{6}$, we get from \eqref{5.10}
\begin{equation*}
\begin{split}
\int_{0}^{T}II_{6}(t)dt
&\leq C\int_{0}^{T}\Big|\int_{\mathbb{T}^{d}}\Cl_{u^{1}}\delta\omega\cdot \varphi dx\Big|dt\\
&\leq C\int_{0}^{T}\|\nabla\delta\omega\|_2\|\varphi\|_2dt\\
&\leq C\|\varphi\|_{L^{2}(0,T;L^{2})}\|\nabla \delta\omega\|_{L^{2}(0,T;L^{2})}\\
&\leq C\|\varphi\|_{L^{4}(0,T;L^{2})}\|\nabla \delta\omega\|_{L^{2}(0,T;L^{2})}\\
&\leq c(T)\|\nabla \delta u\|_{L^{2}(0,T;L^{2})}\|\nabla\delta\omega\|_{L^{2}(0,T;L^{2})}.
\end{split}
\end{equation*}
For $II_{7}$, using \eqref{5.2} and $\|A_{i}\|_{\infty}<\infty(i=1,2)$, we have
\begin{equation*}
\begin{split}
II_{7}(t)&\leq\gamma\Big|\int_{\mathbb{T}^{d}}\Dv\Big(({}^{T}A_{2}A_{2}-{}^{T}A_{1}A_{1})\nabla\bar{\omega}^{2}\Big)\cdot\delta\omega dx\Big|\\
&\leq C \Big|\int_{\mathbb{T}^{d}}\Dv\Big(({}^{T}\delta AA_{2}+{}^{T}A_{1}\delta A)\nabla\bar{\omega}^{2}\Big)
\cdot\delta\omega dx\Big|\\
&\leq C\Big|\int_{\mathbb{T}^{d}}\big|{}^{T}\delta AA_{2}+{}^{T}A_{1}\delta A\big|\big|\nabla\bar{\omega}^{2}\big|
\big|\nabla\delta\omega\big| dx\Big|\\
&\leq C\|t^{-\frac{1}{2}}\delta A\|_{2}\|t^{\frac{1}{2}}\nabla\bar{\omega}^{2}\|_{\infty}
\|\nabla\delta\omega\|_{2},
\end{split}
\end{equation*}
which along with \eqref{5.6} and \eqref{5.13} yields that
\begin{equation*}
\begin{split}
\int_{0}^{T}II_{7}(t)dt
&\leq C\int_{0}^{T}\|t^{-\frac{1}{2}}\delta A\|_{2}\|t^{\frac{1}{2}}\nabla\bar{\omega}^{2}\|_{\infty}
\|\nabla\delta\omega\|_{2}dt\\
&\leq C\|t^{-\frac{1}{2}}\delta A\|_{L^{\infty}(0,T;L^{2})}
\|t^{\frac{1}{2}}\nabla\bar{\omega}^{2}\|_{L^{2}(0,T;L^{\infty})}
\|\nabla\delta\omega\|_{L^{2}(0,T;L^{2})}\\
&\leq c(T)\|\nabla\delta u\|_{L^{2}(0,T;L^{2})}\|\nabla\delta\omega\|_{L^{2}(0,T;L^{2})}.
\end{split}
\end{equation*}
For $II_{8}$, note that $\Dv u=0$, then $J=1$, $a_i=A_i(i=1,2)$.  From Piola identity \eqref{Piola},  we get
\begin{equation}\label{4.43}a_{2j,k}^{k}=a_{1j,k}^{k}=0.\end{equation} Combining  with \eqref{5.2-1}, \eqref{4.43},  H\"{o}lder's inequality and $\|A_{i}\|_{\infty}<\infty(i=1,2)$, we have
\begin{equation*}
\begin{split}
II_{8}(t)&\leq\kappa\Big|\int_{\mathbb{T}^{d}}\Big(\big(\nabla_{u^{2}}\Dv_{u^{2}}-\nabla_{u^{1}}\Dv_{u^{1}}\big)\bar{\omega}^{2}\Big)
\cdot\delta\omega dx\Big|\\
&\leq C\Big|\int_{\mathbb{T}^{d}}\Big({}^{T}A_{2}\nabla\Dv(A_{2}\bar{\omega}^{2})
-{}^{T}A_{1}\nabla\Dv(A_{1}\bar{\omega}^{2})\Big)\cdot\delta\omega dx\Big|\\
&\leq C\Big|\int_{\mathbb{T}^{d}}\Big({}^{T}a_{2}\nabla\Dv(A_{2}\bar{\omega}^{2})
-{}^{T}a_{1}\nabla\Dv(A_{1}\bar{\omega}^{2})\Big)\cdot\delta\omega dx\Big|\\
&\leq C\Big|\int_{\mathbb{T}^{d}}\Big({}^{T}a_{2}\nabla({}^{T}A_{2}:\nabla \bar{\omega}^{2})
-{}^{T}a_{1}\nabla({}^{T}A_{1}:\nabla \bar{\omega}^{2}))\cdot\delta\omega dx\Big|\\
&\leq C\Big|\int_{\mathbb{T}^{d}}\Big({}^{T}a^{j}_{2k}\partial_{k}({}^{T}A_{2}:\nabla \bar{\omega}^{2})-{}^{T}a^{j}_{1k}\partial_{k}(A_{1}:\nabla \bar{\omega}^{2})\Big)\cdot(\delta\omega)^{j} dx\Big|\\
&\leq \Big|\int_{\mathbb{T}^{d}}\Big(a_{2j}^{k}\partial_{k}({}^{T}A_{2}:\nabla \bar{\omega}^{2})-a_{1j}^{k}\partial_{k}(A_{1}:\nabla \bar{\omega}^{2})\Big)\cdot(\delta\omega)^{j} dx\Big|\\
&\leq C\Big|\int_{\mathbb{T}^{d}}\partial_{k}\Big(a_{2j}^{k}({}^{T}A_{2}:\nabla \bar{\omega}^{2})-a_{1j}^{k}(A_{1}:\nabla \bar{\omega}^{2})\Big)\cdot(\delta\omega)^{j} dx\Big|\\
&\leq C\Big|\int_{\mathbb{T}^{d}}\Dv\Big(a_{2}\cdot({}^{T}A_{2}:\nabla \bar{\omega}^{2})
-a_{1}\cdot({}^{T}A_{1}:\nabla \bar{\omega}^{2})\Big)\cdot\delta\omega dx\Big|\\
&\leq C\Big|\kappa\int_{\mathbb{T}^{d}}\Big(a_{2}\cdot({}^{T}A_{2}:\nabla \bar{\omega}^{2})
-a_{1}\cdot({}^{T}A_{1}:\nabla \bar{\omega}^{2})\Big):\nabla\delta\omega dx\Big|\\
&\leq C\Big|\int_{\mathbb{T}^{d}}\Big((a_{2}-a_{1})\cdot({}^{T}A_{2}:\nabla \bar{\omega}^{2})
+a_{1}\cdot({}^{T}\delta A:\nabla \bar{\omega}^{2})\Big):\nabla\delta\omega dx\Big|\\
&\leq C\Big|\int_{\mathbb{T}^{d}}\Big(\delta A\cdot({}^{T}A_{2}:\nabla \bar{\omega}^{2})
+A_{1}\cdot({}^{T}\delta A:\nabla \bar{\omega}^{2})\Big):\nabla\delta\omega dx\Big|\\
&\leq C\Big|\int_{\mathbb{T}^{d}}\Big(\delta A\cdot({}^{T}A_{2}:\nabla \bar{\omega}^{2})
+A_{1}\cdot({}^{T}\delta A:\nabla \bar{\omega}^{2})\Big):\nabla\delta\omega dx\Big|\\
&\leq C\Big(\|t^{-\frac{1}{2}}\delta A\|_{L^{2}}\|{}^{T}A_{2}\|_{L^{\infty}}
\|t^{\frac{1}{2}}\nabla \bar{\omega}^{2}\|_{L^{\infty}}
+\|t^{-\frac{1}{2}}\delta A\|_{L^{2}}\|A_{1}\|_{L^{\infty}}
\|t^{\frac{1}{2}}\nabla \bar{\omega}^{2}\|_{L^{\infty}}\Big)\|\nabla\delta\omega\|_{L^{2}}\\
&\leq C\|t^{-\frac{1}{2}}\delta A\|_{L^{2}}\|t^{\frac{1}{2}}\nabla \bar{\omega}^{2}\|_{L^{\infty}}
\|\nabla\delta\omega\|_{L^{2}},
\end{split}
\end{equation*}
where ${}^{T}a^{j}_{ik}$ denotes the $j^{th}$ row and $k^{th}$ column
component of the matrix ${}^{T}a_{i}(i=1,2)$. From \eqref{5.6} and \eqref{5.13}, we obtain
\begin{equation*}
\begin{split}
\int_{0}^{T}II_{8}(t)dt&\leq  C\|t^{-\frac{1}{2}}\delta A\|_{L^{\infty}(0,T;L^{2})}
\|t^{\frac{1}{2}}\nabla \bar{\omega}^{2}\|_{L^{2}(0,T;L^{\infty})}
\|\nabla\delta\omega\|_{L^{2}(0,T;L^{2})}\\
&\leq c(T)\|\nabla\delta u\|_{L^{2}(0,T;L^{2})}
\|\nabla\delta\omega\|_{L^{2}(0,T;L^{2})}.
\end{split}
\end{equation*}
Finally, for $II_{9}$,  using \eqref{5.2}, \eqref{curl} and  H\"{o}lder's inequality, we have
\begin{equation*}
\begin{split}
II_{9}(t)
&\leq2\chi\Big|\int_{\mathbb{T}^{d}}\Big(\Cl_{u^{2}}-\Cl_{u^{1}}\Big)\bar{u}^{2}\cdot\delta\omega dx\Big|\\
&\leq2\chi\Big|\int_{\mathbb{T}^{d}}\Big(\Cl_{u^{2}}-\Cl_{u^{1}}\Big)\delta\omega\cdot \bar{u}^{2}dx\Big|\\
&\leq C\Big|\int_{\mathbb{T}^{d}}\Big(\varepsilon_{ijk}A_{2j}^{k}(\delta\omega)_{,s}^{k}
-\varepsilon_{ijk}A_{1j}^{k}(\delta\omega)_{,s}^{k}\Big)\cdot(\bar{u}^{2})^{i} dx\Big|\\
&\leq C\Big|\int_{\mathbb{T}^{d}}\varepsilon_{ijk}\delta A_{j}^{k}(\delta\omega)_{,s}^{k}\cdot(\bar{u}^{2})^{i} dx\Big|\\
&\leq C\Big|\int_{\mathbb{T}^{d}}\varepsilon_{ijk}\delta A_{j}^{k}(\delta\omega)_{,s}^{k}\cdot(\bar{u}^{2})^{i} dx\Big|\\
&\leq C\|t^{-\frac{1}{2}}\delta A\|_{L^{2}}
\|t^{\frac{1}{2}}\bar{u}^{2}\|_{L^{\infty}}
\|\nabla\delta\omega\|_{L^{2}},
\end{split}
\end{equation*}
from which,  together with \eqref{5.6} and \eqref{5.13} yields that
\begin{equation*}
\begin{split}
\int_{0}^{T}II_{9}(t)dt
&\leq C\int_{0}^{T}\|t^{-\frac{1}{2}}\delta A\|_{L^{2}}
\|t^{\frac{1}{2}}\bar{u}^{2}\|_{L^{\infty}}
\|\nabla\delta\omega\|_{L^{2}}dt \\
&\leq C\|t^{-\frac{1}{2}}\delta A\|_{L^{\infty}(0,T;L^{2})}
\|t^{\frac{1}{2}}\bar{u}^{2}\|_{L^{2}(0,T;L^{\infty})}
\|\nabla\delta\omega\|_{L^{2}(0,T;L^{2})} \\
&\leq CT^{\frac{3}{4}}\|t^{-\frac{1}{2}}\delta A\|_{L^{\infty}(0,T;L^{2})}
\|\bar{u}^{2}\|_{L^{4}(0,T;L^{\infty})}
\|\nabla\delta\omega\|_{L^{2}(0,T;L^{2})} \\
&\leq c(T)\|\nabla \delta u\|_{L^{2}(0,T;L^{2})}\|\nabla\delta\omega\|_{L^{2}(0,T\times\mathbb{T}^{d})}.
\end{split}
\end{equation*}
So altogether, this gives for all small enough $T>0$,
\begin{align}\label{5.27}\begin{split}
\sup_{t\in[0,T]}\|(\sqrt{\rho_{0}}\phi,\sqrt{\rho_{0}}\delta\omega)\|_{2}^{2}&+\|(\nabla \delta u,\nabla \delta \omega)\|^{2}_{L^{2}(0,T;L^{2})}\leq c(T)\|(\nabla \delta u,\nabla \delta \omega)\|^{2}_{L^{2}(0,T;L^{2})}.
\end{split}\end{align}
Combining with \eqref{5.10}, we conclude that
$$
\|(\nabla \delta u,\nabla \delta \omega)\|^{2}_{L^{2}(0,T;L^{2})}\leq c(T)\|(\nabla \delta u,\nabla \delta \omega)\|^{2}_{L^{2}(0,T;L^{2})}.
$$
Hence $\nabla\delta u=\nabla\delta\omega\equiv0$ on $[0,T]\times\mathbb{T}^{d}$ if $T$ is small enough.
Then, plugging that information into \eqref{5.27} yields
$$\|(\sqrt{\rho_{0}}\phi,\sqrt{\rho_{0}}\delta\omega)\|_{L^{\infty}(0,T;L^{2})}^{2}+\|(\nabla \phi,\nabla\delta\omega)\|_{L^{2}(0,T\times\mathbb{T}^{d})}^{2}=0.$$
Thus, we get $\delta\omega\equiv0$ on $[0,T]\times\mathbb{T}^{d}$ if $T$ is small enough. Combining with Lemma \ref{l:div-dom-1} finally implies that $\phi\equiv0$ on $[0,T]\times\mathbb{T}^{d}$, and \eqref{5.10} clearly yields $\varphi\equiv0$. Therefore,  for small enough $T>0,$ we finally conclude that
$$\bar{u}^{1}=\bar{u}^{2},\quad \bar{\omega}^{1}=\bar{\omega}^{2} \quad on \quad [0,T]\times\mathbb{T}^{d}.$$
Reverting to Eulerian coordinates, we conclude that the two solutions of system \eqref{eq:NA} coincide on $[0,T]\times\mathbb{T}^{d}$. Then standard connectivity arguments yield uniqueness on the whole
$\mathbb{R}^{+}.$

\section{Appendix}
We here list  the useful lemmas and inequalities that have been used several times in the proof of  uniqueness.

\begin{Lemma}\label{l:div}\cite{DM1,DM} Let $A$ be a matrix valued function on $[0,T]\times \mathbb{T}^d$
satisfying
\begin{equation}\label{eq:detA}
\det A\equiv1.
\end{equation}
There exists a  constant $c$ depending only on $d$, such that if
 \begin{equation}\label{eq:smallA}
 \|Id-A\|_{L^\infty(0,T;L^\infty)}+\|A_t\|_{L^2(0,T;L^6)}\leq c,
 \end{equation}
 then for all function $g:[0,T]\times \mathbb{T}^d\to \mathbb{R}$ satisfying
$g\in L^2(0,T\times \mathbb{T}^d)$ and
$$ g=\div R\ \hbox{ with }\ R\in L^4(0,T;L^2)\ \hbox{ and }\ R_t\in L^{4/3}(0,T;L^{3/2}),
$$
the equation
$$\div(Aw)=g \quad\hbox{in}\quad [0,T]\times \mathbb{T}^d$$
admits a solution $w$ in the  space
$$X_T:=\Bigl\{v\in L^4(0,T;L^2(\mathbb{T}^d))\:, \: \nabla v\in L^2(0,T;L^2(\mathbb{T}^d))
\  \hbox{and}\ v_t\in L^{4/3}(0,T;L^{3/2}(\mathbb{T}^d))\Bigr\}
$$
satisfying the following  inequalities for some  constant $C=C(d)$:
\begin{equation}\label{eq:diveq}\begin{array}{c}
\|w\|_{L^4(0,T;L^2)}\leq C\|R\|_{L^4(0,T;L^2)}, \quad
\|\nabla w\|_{L^2(0,T;L^2)}\leq C\|g\|_{L^2(0,T;L^2)}\\[1ex] \hbox{\rm and }\
\|w_t\|_{L^{4/3}(0,T;L^{3/2})}\leq C\|R\|_{L^4(0,T;L^{2})}+C\|R_t\|_{L^{4/3}(0,T;L^{3/2})}.\end{array}\end{equation}
\end{Lemma}
In the bounded domain case, the previous lemma can be adapted as follows.
\begin{Lemma}\label{l:div-dom}\cite{DM1,DM} Let $\Omega$ be a $\cC^2$ bounded domain of $\mathbb{R}^d,$ and  $A,$  a matrix valued function on $[0,T]\times\Omega$
satisfying  \eqref{eq:detA}. If  \eqref{eq:smallA} is fulfilled
 then for all function $R:[0,T]\times\Omega\to \mathbb{R}^d$ satisfying
$\div R\in L^2(0,T\times\Omega),$ $R\in L^4(0,T;L^2),$ $R_t\in L^{4/3}(0,T;L^{3/2})$ and $R\cdot n\equiv0$ on $(0,T)\times\d\Omega,$
the equation
$$\div(Aw)=\div R=:g \quad\hbox{in}\quad [0,T]\times \Omega$$
admits a solution in the  space
$$X_T:=\Bigl\{v\in L^2(0,T;H^1_0(\Omega))\:, \:  v\in L^4(0,T;L^2(\Omega))
\  \hbox{and}\ v_t\in L^{4/3}(0,T;L^{3/2}(\Omega))\Bigr\},
$$
that satisfies Inequalities \eqref{eq:diveq}.
\end{Lemma}
\begin{Lemma}\label{l:div-dom-1}\cite{BD} Let $a :(0,1)^{d}\rightarrow \mathbb{R}$ be a nonnegative and nonzero measurable function. Then we have for all $z$ in $H^{1}(\mathbb{T}^{d})$,
\begin{equation*}
\|z\|_{2}\leq\frac{1}{M}\Big|\int_{\mathbb{T}^{d}}azdx\Big|+\Big(1+\frac{1}{M}\|M-a\|_{2}\Big)\|\nabla z\|_{2} \quad \text{with} \quad M:=\int_{\mathbb{T}^{d}}adx.
\end{equation*}
Furthermore, in dimension $d=2,$ there exists an absolute constant $C$ so that
\begin{equation*}
\|z\|_{2}\leq\frac{1}{M}\Big|\int_{\mathbb{T}^{d}}azdx\Big|+C\log^{\frac{1}{2}}\Big(e
+\frac{\|M-a\|_{2}}{M}\Big)\|\nabla z\|_{2}.
\end{equation*}
\end{Lemma}
\begin{center}

\end{center}
\end{document}